\documentclass[11pt]{amsart}
\numberwithin{equation}{section}
\usepackage[dvips]{graphicx}
\usepackage{amssymb}
\usepackage{amsmath}
\usepackage{latexsym}
\newcommand{\eps}{\varepsilon}

\newcommand{\N}{{\mathbb N}}
\newcommand{\C}{{\mathbb C}}
\newcommand{\Z}{{\mathbb Z}}
\newcommand{\R}{{\mathbb R}}
\newcommand{\J}{{\mathcal J}}

\newcommand{\bwd}{Baker wandering domain}
\newcommand{\tef}{transcendental entire function}
\newcommand{\tmf}{transcendental meromorphic function}
\newcommand{\bcc}{bounded complementary component}

\newcommand{\wt}{\widetilde}

\def\R{\mathbb R}
\def\C{\mathbb C}
\def\CC{\widehat{\mathbb C}}
\def\N{\mathbb N}

\def\OO{{O}}
\def\im{\operatorname{Im}}
\def\re{\operatorname{Re}}

\def\arg{\operatorname{arg}}

\def\meas{\operatorname{meas}}
\def\supp{\operatorname{supp}}
\def\interior{\operatorname{int}}
\newtheorem{la}{Lemma}[section]
\newtheorem{thm}{Theorem}[section]
\newtheorem{prop}{Proposition}[section]
\theoremstyle{definition}
\newtheorem{defin}{Definition}[section]
\theoremstyle{remark}
\theoremstyle{remark}
\newtheorem*{rem}{Remark}
\newtheorem{remark}{Remark}
\newtheorem{ex}{Example}[section]
\newtheorem{ack}{Acknowledgement}
\title[Direct and logarithmic singularities]{Dynamics 
of meromorphic functions with direct 
or logarithmic singularities}
\thanks{The authors were 
supported by a London Mathematical Society
Scheme~1 grant.
The first author was also supported by the G.I.F.,
the German--Israeli Foundation for Scientific Research and
Development, Grant G-809-234.6/2003, 
the EU Research Training Network CODY
and the ESF Research Networking Programme HCAA.}
\author{Walter Bergweiler}
\address{Mathematisches Seminar,
Christian--Albrechts--Universit\"at zu Kiel,
Lude\-wig--Meyn--Str.~4,
D--24098 Kiel,
Germany}
\email{bergweiler@math.uni-kiel.de}
\author{Philip J.\ Rippon}
\address{Department of Mathematics, The Open University, Walton Hall, 
Milton Keynes MK7 6AA, United Kingdom}
\email{p.j.rippon@open.ac.uk}
\author{Gwyneth M.\ Stallard }
\address{Department of Mathematics, The Open University, Walton Hall, 
Milton Keynes MK7 6AA, United Kingdom}
\email{g.m.stallard@open.ac.uk}
\date{\today}
\begin{document}
\begin{abstract}
Let $f$ be a transcendental meromorphic function
and denote by $J(f)$ the Julia set and by $I(f)$ the escaping set.
We show that if $f$ has a direct singularity over infinity, then
$I(f)$ has an unbounded component and  $I(f)\cap J(f)$ contains 
continua. Moreover, under this hypothesis
$I(f)\cap J(f)$ has an unbounded component if
and only if $f$ has no Baker wandering domain.
If $f$ has a logarithmic singularity over infinity, then
the upper box dimension of $I(f)\cap J(f)$ is $2$ and the Hausdorff
dimension of $J(f)$ is strictly greater than~$1$.
The above theorems are deduced from more general
results concerning functions which have ``direct or logarithmic tracts'',
but which need not be meromorphic in the plane.
These results are obtained by using a generalization of Wiman-Valiron theory.
This method is also applied to complex differential equations.
\end{abstract}
\maketitle
\section{Introduction}
For a function $f$ meromorphic in the plane 
the {\em Fatou set} $F(f)$ is defined as the set where the iterates $f^n$
of $f$ are defined and form a normal family, and
the  {\em Julia set} $J(f)$ is its complement.
The {\em escaping set} $I(f)$ is defined as the set of all $z\in\C$ for which
$(f^n(z))$ is defined and $f^n(z)\to\infty$ as $n\to\infty$.

While the main objects studied in complex dynamics are the Fatou
and Julia sets,
the escaping set also plays a major role in the iteration theory of
entire functions, beginning with a paper by Eremenko~\cite{Ere89} 
who proved that if $f$ is an entire transcendental function, then
$I(f)\neq \emptyset$, $\partial I(f)=J(f)$, $I(f)\cap J(f)
\neq \emptyset$ and 
$\overline{I(f)}$ has no bounded components.
Eremenko conjectured that in fact all components of $I(f)$ 
are unbounded. While this conjecture 
is still open, it is known that $I(f)$ has at least
one unbounded component~\cite{RS05}.
A detailed study of the escaping set of an exponential function was given
by Schleicher and Zimmer~\cite{SZ}.
In particular, they confirmed Eremenko's conjecture for such functions.
Since then the conjecture has also been proved
for certain more general classes of entire functions
by Bara\'nski~\cite{Bar},
Rempe~\cite{Rem1}
and Rottenfu{\ss}er, R\"uckert, 
Rempe and Schleicher~\cite{RRRS}.
Actually it was shown in~\cite{Bar,RRRS,SZ} that 
-- in the classes of functions considered --
every point  in $I(f)$ can be connected
to infinity by a curve in $I(f)$.
However, in~\cite{RRRS}, an example of an entire function 
was given for which this is not
true, even though, by the results in~\cite{Rem1},
Eremenko's conjecture is satisfied for
this function. This disproves a stronger form of Eremenko's conjecture.
Among further
results related to the escaping set of an entire
function we mention~\cite{Kar,McM,Rem,Schl}.

For a meromorphic function $f$ the set
$I(f)$ was first considered by Dom\'{\i}nguez~\cite{Dom98}.
She proved that again $I(f)\neq \emptyset$, $\partial I(f)=J(f)$
and $I(f)\cap J(f) \neq \emptyset$.
On the other hand, in this situation
$\overline{I(f)}$ need not have unbounded components. For example, 
for $f(z)=\frac12 \tan z$ 
we have $\overline{I(f)}=J(f)$, and this set is totally disconnected.

In general it can be said that the set $I(f)$ is much less useful
for meromorphic functions with infinitely many poles than for 
entire functions.
For example, in many of the recent studies
(e.g.,~\cite{RS00,RS05})
of the set $I(f)$
for entire $f$ it turned out to be useful to consider the 
rate of escape to infinity for points in $I(f)$. 
For functions with infinitely many poles such escape rates are
of limited use, since a point can escape to infinity 
arbitrarily fast by ``jumping from pole to pole''.
However, 
some of the results about entire functions have been carried over
to meromorphic functions with finitely many 
poles~\cite{Dom98,RS00,RS05b,RS06}.

Here we single out a large class of meromorphic functions which may 
have infinitely many poles, but which still have dynamical properties 
similar to entire functions.
More specifically, we consider meromorphic functions which
have a direct singularity over infinity. (The definition 
of a direct singularity will be given at the beginning of Section~\ref{direct}.)
\begin{thm} \label{unboundedcomp}
Let $f$ be a meromorphic function with a direct singularity over infinity.
Then $I(f)$ has an unbounded component.
\end{thm}
\begin{thm} \label{continua}
Let $f$ be a meromorphic function with a direct singularity over infinity.
Then $I(f)\cap J(f)$ contains continua.
\end{thm}
In general, $I(f)\cap J(f)$ need not have unbounded components, but 
the case where this does not happen can be characterized. 
Recall that if $f$ is a {\tmf} and if $U$ is a component of $F(f)$,
then there exists, for each $n\in \N$, a component of $F(f)$, which
we call $U_n$, such that $f^n(U) \subset U_n$.
If $U_m\neq U_n$ whenever $m\neq n$, then $U$ is called a 
{\it wandering domain}.
We use the name
{\it Baker wandering domain} to denote a wandering domain $U$ of $F(f)$
such that, for $n$ large enough, $U_n$ is a bounded multiply connected
component of $F(f)$ which surrounds~0, and $U_n \to \infty$ as $n\to\infty$.
The first example of such a wandering
domain was given by Baker~\cite{Bak76}.
\begin{thm} \label{unboundedifnobwd}
Let $f$ be a meromorphic function with a direct singularity over infinity.
Then $I(f)\cap J(f)$ has an unbounded component if and only if 
$f$ has no Baker wandering domain.
\end{thm}

The proofs of the above results for entire functions 
(and also the proof that $I(f)\neq \emptyset$ for entire $f$) use 
Wiman-Valiron theory or the maximum principle and thus
do not carry over to meromorphic functions with poles.
Our main tools are results whose statements 
are similar to those of Wiman-Valiron theory~\cite{Hay74},
or of Macintyre's theory of flat regions~\cite{Mac38}.
We will develop these results in Section~\ref{direct},
but the methods are
quite different from those of Wiman, Valiron and Macintyre.

An important special class of meromorphic functions is the 
Eremenko-Lyubich
class $B$ consisting of those functions for which there exists
$R>0$ such that the inverse function of $f$ has no singularity
over $\{z\in \C:|z|>R\}$.
The main tool used to study this class of functions is the logarithmic
change of variable, introduced  to the
subject by Eremenko and Lyubich~\cite{EL}.
In particular, they proved that if a transcendental entire function 
$f$ belongs to the class $B$, then $I(f)\subset J(f)$.

Entire functions in the class $B$ have
a logarithmic singularity over infinity.
Many results proved about the class $B$ 
hold more generally for functions with a 
logarithmic singularity over infinity, with 
proofs carrying over to this more general setting without change.
This observation 
is not new.
For example, the results in~\cite{Rem,RRRS} are already
stated for functions with a logarithmic singularity over infinity,
and not just for entire functions in class~$B$.

However, there are some results for class $B$
where the generalization from
``entire in class $B$'' to 
``meromorphic with a logarithmic singularity'' 
is not so obvious.
In particular, this applies to the following result, 
where the proof for entire functions used Wiman-Valiron theory.
\begin{thm} \label{box}
Let $f$ be a meromorphic function with a logarithmic
singularity over infinity.
Then the upper box dimension of $I(f)\cap J(f)$ is $2$ and
the packing dimension of $J(f)$ is~$2$.
\end{thm}
Theorem~\ref{box} generalizes~\cite[Theorem~1.1]{RS05c} and we briefly outline
here the changes in the proof needed in the present context.
Similarly we sketch the changes needed in the proof of~\cite[Theorem~3]{RS06}
to give the following result.
\begin{thm} \label{hausdorff}
Let $f$ be a meromorphic function with a logarithmic
singularity over infinity.
Then the Hausdorff dimension of $J(f)$ is 
strictly greater than~$1$.
\end{thm}
The paper is organized as follows.
In Section~\ref{direct} we recall the classification of singularities of the inverse
function and state the results of Wiman-Valiron type which will be 
used when dealing with direct singularities over infinity.
The proof of these results will be deferred until
Sections~\mbox{\ref{proofgrowth}--\ref{proof1}}.
In Section~\ref{iterationtract} we discuss iteration within
a ``tract'' associated to the direct singularity over infinity
and prove two results about fast escaping points whose orbits 
lie eventually in the tract (Theorems~\ref{th3} 
and~\ref{th4}) from which Theorem~\ref{unboundedcomp}
immediately follows.
While the results of Section~\ref{iterationtract} only require that the 
function is defined in a ``tract'', 
the following two sections
are devoted to functions meromorphic in the plane.
In particular, Theorems~\ref{continua} and~\ref{unboundedifnobwd}
will be corollaries to two results 
(Theorems~\ref{bwd} and~\ref{nobwd})
proved in Section~\ref{functionwithbwd}.
In Section~\ref{classB} we consider functions with a logarithmic
tract and prove results from which Theorems~\ref{box} and~\ref{hausdorff}
follow.
In Section~\ref{teichm} we discuss results of Teichm\"uller and
Selberg which show that, under suitable additional hypotheses,
a meromorphic function has a logarithmic or direct singularity if it has
few poles in the sense of Nevanlinna theory.
Some examples will be discussed in Section~\ref{examples}.

Finally we note that the results of
Section~\ref{direct}  may be applied not only in complex dynamics, but also to
complex differential equations. One such application will be given
in Section~\ref{diffeq}.
\begin{rem}
After completion of the paper, Alex Eremenko kindly informed us that the 
idea of using Wiman-Valiron theory for meromorphic
functions with direct tracts appears already
in his paper~\cite{Ere82} 
dealing with differential equations of Briot-Bouquet type.
This paper contains a statement equivalent to our
Theorem~\ref{th1}, but according to A. Eremenko his proof
of this statement contains a gap. Our Theorem~\ref{th1} fills this
gap.
\end{rem}
\begin{ack}
We thank Alex Eremenko, Katsuya Ishizaki, Jim Langley,
Lasse Rempe and the referee for very valuable comments.
\end{ack}

\section{Direct and logarithmic tracts}
\label{direct}
We denote the open disc of radius $r$ around a point $a\in\C$ by
$D(a,r)$ and the closed disc by $\overline{D}(a,r)$.
The open disc of radius $r$ around a point $a\in\CC=\C\cup\{\infty\}$ with
respect to the spherical metric is denoted by $D_\chi(a,r)$.

We recall the classification of singularities of the inverse
function of a meromorphic function due to Iversen~\cite{Ive};
see~\cite[p.~289]{Nev53}.
Let $f$ be meromorphic in the plane and let $a\in \CC$.
For $r>0$  let $U_r$ be a component of the preimage
$f^{-1}\left(D_\chi(a,r)\right)$, chosen in such a way
that $r_1<r_2$ implies $U_{r_1}\subset U_{r_2}$. 
Then there are two possibilities:
\begin{itemize}
\item[(a)] $\bigcap_{r>0}U_r=\{z\}$ for some $z\in \C$, or
\item[(b)] $\bigcap_{r>0}U_r=\emptyset.$
\end{itemize}
In case (a) we have $a=f(z)$.
If $a\in \C$ and $f'(z)\neq 0$
or if $a=\infty$ and $z$ is a simple pole of $f$,
then $z$ is called an {\em ordinary point}.
If $a\in \C$ and $f'(z)= 0$
or if $a=\infty$ and $z$ is a multiple pole of $f$,
then $z$ is called
a {\em critical point} and $a$ is called a {\em critical value}.
We also say $f^{-1}$ has an {\em algebraic singularity} over~$a$.

In case (b) we say that our choice $r\mapsto U_r$ defines
a {\em transcendental singularity} of $f^{-1}$ over~$a$.
A transcendental singularity over $a$ is called {\em direct} if for
some $r>0$ we have $f(z)\neq a$ for $z\in U_r$. Otherwise it is called
{\em indirect}. A direct singularity is called {\em logarithmic} if
the restriction $f:U_r\to D_\chi(a,r)\backslash\{ a\}$ is a 
universal covering for some $r>0$.

Note that in case (b) there exists a curve $\gamma$ tending to $\infty$
such that $f(z)$ tends to $a$ as $z$ tends to infinity along $\gamma$.
A value $a$ for which such a curve exists is called an
{\em asymptotic value}. 
In turn, if $a$ is an asymptotic value and $\gamma$ is the corresponding
curve, then by choosing 
$U_r$ as the component of 
$f^{-1}\left(D_\chi(a,r)\right)$ which contains 
the ``tail'' of this curve we 
obtain a transcendental singularity over~$a$.
Thus we see that $f$ has a transcendental singularity over $a$
if and only if $a$ is an asymptotic value.

We note that over a point $a$ there can be several singularities,
of the same or of
different types, and there may also be ordinary points over~$a$.
For example, the function $f(z)=\Gamma(z)\exp(z^2)$ has 
two transcendental 
singularities over $\infty$, with corresponding asymptotic paths being
the positive and negative real axes. The singularity corresponding
to the positive real axis is direct while the one corresponding
to the negative real axis is indirect.

The domains $U_r$ appearing in the definition of a direct or logarithmic
singularity are called {\em tracts}. In some of our results we will work
only with the restriction of $f$ to a tract, and it will be irrelevant
how $f$ behaves outside the tract. This motivates the following definition.

\begin{defin} \label{defin1}
Let $D$ be an unbounded domain in $\C$ whose boundary consists
of piecewise smooth curves. Suppose that 
the complement of $D$ is unbounded.
Let $f$ be a complex-valued function whose domain of definition
contains the closure $\overline{D}$ of~$D$. Then $D$ is called a 
{\em direct tract} of $f$ if $f$ is holomorphic in $D$ and
continuous in $\overline{D}$ and if there exists $R>0$ such that
$|f(z)|=R$ for $z\in\partial D$ while 
$|f(z)|>R$ for $z\in D$.
If, in addition, 
the restriction $f:D\to \{z\in\C:|z|>R\}$ is a universal covering,
then $D$ is called a {\em logarithmic tract} of~$f$.
\end{defin}
Note that if $f$ is meromorphic in the plane with a direct singularity
over $\infty$, then $f$ has a direct tract.
In particular, every transcendental entire function has a direct tract.
Similarly, if $f$ has a logarithmic singularity
over $\infty$, then $f$ has a logarithmic tract.

We remark that the Denjoy-Carleman-Ahlfors Theorem~\cite[Section XI.4]{Nev53}
says that a meromorphic function of finite order $\rho$ can have
at most $\max\{1,2\rho\}$ direct singularities, and thus 
at most $\max\{1,2\rho\}$ direct tracts, for each fixed value of~$R$.
In particular, a meromorphic function of order less than $1$ 
can have at most one direct tract.  
Also, this holds 
if $f$ is meromorphic and there is a sequence of Jordan curves $\gamma_n$ 
surrounding the origin and tending to $\infty$ on which $f$
tends to~$\infty$.  In particular, this is the case for a meromorphic
function with a Baker wandering domain.

For a nonconstant subharmonic function  $v:\C\to [0,\infty)$
the function
\begin{equation}\label{defBrv}
B(r,v)=\max_{|z|=r}v(z)
\end{equation}
is increasing, convex in $\log r$ and tends to $\infty$ as $r$
tends to $\infty$; see~\cite[Chapter 2]{HK1}.
Hence
\begin{equation}\label{defarv}
a(r,v)=\frac{d B(r,v)}{d\log r}=rB'(r,v)
\end{equation}
exists except perhaps for a countable set of $r$-values,
and $a(r,v)$ is nondecreasing.

Note that if $f,D,R$ are as in the above definition, then 
the function $v:\C\to [0,\infty)$ defined by 
\begin{equation}\label{definitionv}
v(z)=
\begin{cases} 
\displaystyle \log \frac{|f(z)|}{R} & \text{if $z\in D$},\\[3mm]
0 & \text{if $z\notin D$},
\end{cases}
\end{equation}
is subharmonic.
While for a nonconstant function $v$ subharmonic in the plane 
we only know~\cite[Theorem~2.14]{Hay89} that
$$
\lim_{r\to\infty} \frac{B(r,v)}{\log r} >0,
$$
functions of the form~(\ref{definitionv}) have faster growth.
\begin{thm} \label{growth}
Let $D$ be a direct tract of $f$ and let $v$ be defined
by \textup{(\ref{definitionv})}. Then
\begin{equation} \label{growthB}
\lim_{r\to\infty} \frac{B(r,v)}{\log r} =\infty
\end{equation}
and
\begin{equation} \label{growtha}
\lim_{r\to\infty} a(r,v)=\infty.
\end{equation}
\end{thm} 
This result is due to Fuchs~\cite{Fuc81}.
He proved only~\eqref{growthB}, but~\eqref{growtha} follows 
from~\eqref{growthB} since 
\begin{equation} \label{comparisonBa}
a(r,v)\geq \frac{B(r,v)-B(r_0,v)}{\log (r/r_0)}
\end{equation}
for $r>r_0>0$ by \eqref{defarv}.
We include a proof of Theorem~\ref{growth} in Section~\ref{proofgrowth}
for completeness and because we use the
techniques in the proof later in the paper.

We mention that $a(r,v)$ cannot only be estimated in terms of
$B(r,v)$ from below as in~\eqref{comparisonBa}, but also 
from above; see~Lemma~\ref{growthlemma5}.

Our main tool to study functions with a direct tract is
the following result, which is the main result in the
paper. Here we say that a set $F\subset [1,\infty)$ 
has {\em finite logarithmic measure} if $\int_F dt/t <\infty$.
We also put
\begin{equation} \label{defM}
M_D(r)=\max_{|z|=r, z\in D}|f(z)|=\exp B(r,v).
\end{equation}
\begin{thm} \label{th1}
Let $D$ be a direct tract of $f$ and let $\tau>\frac12$.
Let $v$ be defined by~\textup{(\ref{definitionv})} and
let $z_r$ be a point satisfying $|z_r|=r$ and $v(z_r)=B(r,v)$.
Then there exists a set $F\subset [1,\infty)$ of finite logarithmic
measure such that if $r\in [1,\infty)\setminus F$, then 
$D(z_r,r/a(r,v)^{\tau})\subset D$. Moreover,
\begin{equation}\label{wvlike}
f(z)\sim \left(\frac{z}{z_r}\right)^{a(r,v)}f(z_r)
\quad \text{for} \quad
z\in D\left(z_r,\frac{r}{a(r,v)^\tau}\right)
\end{equation}
and
\begin{equation}\label{wvlike4}
|f(z)| \sim M_D(|z|)
\quad \text{for} \quad
z\in D\left(z_r,\frac{r}{a(r,v)^\tau}\right)
\end{equation}
as $r\to\infty$, $r\notin F$.
\end{thm}
We note that it follows 
from~(\ref{wvlike}) that if $k\in\N$, then
\begin{equation}\label{wvlike3}
f^{(k)}(z)\sim 
\left(\frac{a(r,v)}{z}\right)^k
\left(\frac{z}{z_r}\right)^{a(r,v)}f(z_r)
\quad
\text{for}
\quad
z\in D\left(z_r,\frac{r}{a(r,v)^\tau}\right)
\end{equation}
as $r\to\infty$, $r\notin F$.
We will need~\eqref{wvlike3} in Section~\ref{classB} and Section~\ref{diffeq}.

Since $a(r,v)\to\infty$ as $r\to\infty$
by Theorem~\ref{growth}, since $e^h-1\sim h$ as $h\to 0$ and 
since $\tau>\frac12$ is arbitrary, 
we see that conclusion~(\ref{wvlike}) in Theorem~\ref{th1}
can be replaced by
\begin{equation}\label{wvlike2}
f(z_re^h)\sim e^{a(r,v)h}f(z_r)
\quad
\text{for}
\quad
|h|\leq a(r,v)^{-\tau},
\end{equation}
again as  $r\to\infty$, $r\notin F$.

The asymptotic relations~(\ref{wvlike}), (\ref{wvlike3})
and~(\ref{wvlike2}) are very similar to the main
results of Wiman-Valiron theory (see, e.g.,~\cite{Hay74}),
except that the central index 
is replaced by $a(r,v)$. 
However, our results do not even require 
$f$ to be defined outside of the closure of the tract $D$
whereas the Wiman-Valiron method requires that $f$ is entire,
since it is based on the Taylor series expansion.
It follows from the example in~\cite[pp.~345--346]{Hay74} that we cannot
take $\tau=\frac12$ in Theorem~\ref{th1}.

It is not difficult to see that
$z_r f'(z_r)/f(z_r)$ lies between the left and right derivative
of $B(r,v)$ with respect to $\log r$ and thus 
$$a(r,v)=\frac{z_r f'(z_r)}{f(z_r)}$$
if $a(r,v)$ exists.
(Actually,
by a result of Blumenthal (see~\cite[Section~II.3]{Val23}), 
the set of $r$-values where $a(r,v)$ does not exist is discrete 
if $v=\log |f|$ for some entire function~$f$, and Blumenthal's result
extends to the case where $v$ is as in Theorem~\ref{th1}.)
With the above expression for $a(r,v)$ the relation~(\ref{wvlike2})
was proved by Macintyre~\cite{Mac38} for entire functions.
While it is possible to relax the assumption that $f$ is entire,
it seems essential for his method
that a certain disc around $z_r$  is in the domain of definition 
of~$f$.
The key conclusion of Theorem~\ref{th1} therefore is that 
$D(z_r,r/a(r,v)^{\tau})\subset D$ and most of the proof of 
Theorem~\ref{th1} is devoted to this fact.

As in Wiman-Valiron theory, an important consequence of 
(\ref{wvlike}) or~(\ref{wvlike2}) 
is the following result, which can be deduced 
from them for example by Rouch\'e's theorem.
\begin{thm} \label{th2}
For each $\beta>1$ there exists $\alpha>0$ such that if
$f$, $D$, $v$, $z_r$ and $F$
are as in Theorem~\ref{th1} and if $r\notin F$ is sufficiently large,
then 
$$\left\{z\in\C:\frac{|f(z_r)|}{\beta}\leq |z|\leq \beta |f(z_r)| 
\right\}\subset
f\left(D\left(z_r,\frac{\alpha r}{a(r,v)}\right)\right).$$
More precisely, 
$\log f$ is univalent in $D(z_r,\alpha r/a(r,v))$
and for $\gamma>\pi$ the constant $\alpha$ can be chosen
such that if $r\notin F$ is sufficiently large,
then $\log f \left(D(z_r,\alpha r/a(r,v))\right)$ contains the
rectangle
$$\left\{ z\in\C: \left| \re z -\log |f(z_r)|\right|\leq \log\beta,
\left| \im z -\arg f(z_r)\right| \leq \gamma \right\},
$$
where the branches of $\log$ and $\arg$ are chosen such that
$\im(\log f(z_r))=\arg f(z_r)$.
\end{thm}
Theorem~\ref{th2} is sufficient for our purposes, but we note that this
statement can be strengthened by working with the disc
$D(z_r,r/a(r,v)^{\tau})$, for $\frac12<\tau<1$. 
\section{Iteration in a tract}
\label{iterationtract}
Wiman-Valiron theory was the main tool in Eremenko's proof~\cite{Ere89} 
that $I(f)\neq\emptyset$ if $f$ is a transcendental entire function.
Using 
Theorem~\ref{th2} instead of the Wiman-Valiron method
in his argument, we obtain the following result.
\begin{thm} \label{thnew}
Let $D$ be a direct tract of~$f$. Then there exists $z_0\in D$ such that
$f^n(z_0)\in D$ for all $n\in\N$ and $f^n(z_0)\to\infty$ as $n\to\infty$.
\end{thm}
We emphasize that we are not assuming here that $f$ is defined
in the whole plane, but only that the domain of definition of $f$ 
contains~$\overline{D}$.
However, the result appears to be new even for entire~$f$.
Neither Eremenko's argument based on the Wiman-Valiron method
nor Dom\'{\i}nguez's argument~\cite{Dom98} based on the 
maximum modulus and a theorem of Bohr seem to give, for entire~$f$,
the existence of a point $z_0\in I(f)$ such that $f^n(z_0)$ is in
the same direct tract for all~$n$.
The conclusion of Theorem~\ref{thnew} is known, however,
if $D$ is a logarithmic tract~\cite{Rem2}.
\begin{proof}[Proof of Theorem~\ref{thnew}]
We adapt the argument of Eremenko.
Let $F$ be the exceptional set arising in Theorems~\ref{th1} and~\ref{th2},
with $\beta=4$ in Theorem~\ref{th2}, and $\alpha$ chosen accordingly.
For large $r$ we have
$\left[\frac12 |f(z_r)|,2|f(z_r)|\right]\setminus F\neq\emptyset$
and $ |f(z_r)|>4 r$, by~(\ref{growthB}).
Choosing $r_0\notin F$ large we can thus 
inductively define an increasing sequence $(r_n)$ 
satisfying
$$r_{n+1}\in \left[\tfrac12 |f(z_{r_n})|,2|f(z_{r_n})|\right]\setminus F$$
and $r_n\to\infty$.
We put $\overline{D}_n=
\overline{D}\left(z_{r_{n}},\alpha r_{n}/a(r_{n},v)\right)$.
Choosing $r_0$ large we have
$\overline{D}_{n+1}\subset
\left\{z\in\C:\frac14 |f(z_{r_n})|\leq |z|\leq 4 |f(z_{r_n})| \right\}$
and thus, by Theorem~\ref{th2},
$\overline{D}_{n+1}\subset f\left(\overline{D}_n\right)$ 
for all $n\geq 0$.
Inductively we see that there exists a closed set $C_n\subset \overline{D}_0$
such that 
$f^j(C_n)\subset \overline{D}_j$
for $0\leq j\leq n$ and 
$C_{n+1}\subset C_n$ for all $n\geq 0$.
Choosing $z_0\in\bigcap_{n=1}^\infty C_n$ we see  that
$f^n(z_0)\in D$ for all $n\in\N$ and
$f^n(z_0)\to\infty$ as $n\to\infty$.
\end{proof}
Using~\eqref{wvlike4} we see that the chosen point $z_0$ actually satisfies
\begin{equation} \label{eremenko}
|f^{n+1}(z_0)|\sim M_D(|f^n(z_0)|)
\end{equation}
as $n\to\infty$.

It follows from~(\ref{growthB}) that
$M_D(\rho)> \rho$ for large $\rho$, say $\rho> \rho_0>R$.
Hence $M_D^n(\rho)\to\infty$ as $n\to\infty$ 
for $\rho> \rho_0$.
For such $\rho$ we define
$$
A(f,D,\rho)=\{z\in D: f^n(z)\in D \text{ and } |f^n(z)|
\geq M_D^n(\rho) \text{ for
all } n\in\N \}.
$$
In contrast to similar definitions in~\cite{BH,RS05},
the set $A(f,D,\rho)$ depends on $\rho$.
Note that $f^n(z)\to\infty$ as $n\to\infty$
for $z\in A(f,D,\rho)$.
In particular, if $f$ is meromorphic in $\C$, then
$A(f,D,\rho)\subset I(f)$. 

The following result, proved below, strengthens Theorem~\ref{thnew}.
\begin{thm} \label{th3}
Let $D$ be a direct tract of~$f$. Then $A(f,D,\rho)\neq\emptyset$.
\end{thm}
While the results in~\cite{BH,RS05} yield that for an entire function
$f$ there exists $z_0\in\C$ satisfying 
$|f^{n}(z_0)|\geq  M^n(\rho,f)$, where $M(r,f)=\max_{|z|=r}|f(z)|$,
with the exponent $n$ indicating iteration with respect to the
first variable,
the result that $z_0$ can be chosen such that $f^n(z_0)$
is in the same direct tract of $f$ for every $n\in \N$
is again new even for entire~$f$.
\begin{thm} \label{th4}
Let $D$ be a direct tract of~$f$. Then all components of
$A(f,D,\rho)$  are unbounded.
\end{thm}
\begin{proof}[Proof of Theorem~\ref{th3}]
We follow the proof of Lemma~2 in~\cite{BH}. 
Fix $\rho>\rho_0$.
Since $\log M_D(r)=B(r,v)$ is convex in $\log r$ we deduce 
from~(\ref{growthB}) that 
\begin{equation} \label{g}
M_D(2r) \geq 4 M_D(r)
\end{equation}
for large $r$, say $r\geq r_0>\rho$. 
By Eremenko's argument there exists $z_0\in D$ satisfying~(\ref{eremenko})
such that $f^n(z_0)\in D$ for all $n\in\N$.
We may assume that
\begin{equation} \label{h}
|f^{n+1}(z_0)|\geq \tfrac{1}{2} M_D\left(|f^n(z_0)|\right)
\end{equation}
for all $n\in \N$ and that
\begin{equation} \label{i}
|z_0|\geq 2 r_0,
\end{equation}
because otherwise we can replace $z_0$ by
$f^k(z_0)$ for a sufficiently large~$k$.
We shall prove by induction that 
\begin{equation} \label{induction}
|f^{n}(z_0)| \geq 2M_D^n(\rho)
\end{equation}
for all $n\geq 0$.
Because of~(\ref{i}) and since $r_0\geq \rho$ we see that
(\ref{induction}) holds for $n=0$. 
Suppose it holds for some $n\geq 0$.
Combining this with~(\ref{g}) and~(\ref{h}) we deduce that
\[
|f^{n+1}(z_0)|\geq \tfrac{1}{2} M_D\left(|f^n(z_0)|\right)
\geq 2 M_D\left(\tfrac12 |f^n(z_0)|\right)
\geq 2 M_D\left(M_D^n(\rho)\right)
=2 M_D^{n+1}(\rho),
\]
and thus~(\ref{induction}) also holds with $n$ replaced by $n+1$.
Of course, it follows from~(\ref{induction}) that $z_0\in
A(f,D,\rho)$.
\end{proof}
\begin{proof}[Proof of Theorem~\ref{th4}]
Let $z_0\in A(f,D,\rho)$. We follow the argument in~\cite{RS05}
and denote by $L_n$ the component of 
$f^{-n}(\C\setminus D(0,M_D^n(\rho)))$ that contains $z_0$. 
Then $L_n$ is closed and unbounded.
In fact, since 
$|f(z)|=R<\rho\leq M_D^k(\rho)$ for $z\in \partial D$ and $0\leq k<n$, we 
have $f^k(z)\in D$ for $z\in L_n$ and $0\leq k<n$; 
in particular, $L_n\subset D$.
We claim that $L_{n+1}\subset L_n$.
Otherwise there exists  $z'\in L_{n+1}$ 
with $f^n(z')\in D(0,M_D^n(\rho))$. This implies that
$|f^{n+1}(z')|< M_D(M_D^n(\rho)))=M_D^{n+1}(\rho)$ 
so that $f^{n+1}(z')\notin \C\setminus D(0,M_D^{n+1}(\rho))$,
contradicting $z'\in L_{n+1}$.
Hence $L_{n+1}\subset L_n$.
As in~\cite{RS05} we conclude that
$$K=\bigcap_{n=1}^\infty L_n\cup\{\infty\}$$
is a closed connected subset of $\CC$ containing $z_0$ and $\infty$.
It follows from the
construction that we also have $K\setminus\{\infty\}\subset A(f,D,\rho)$.
The component of $A(f,D,\rho)$ containing $z_0$ thus also 
contains a component of $K\setminus\{\infty\}$, and hence it is unbounded;
see~\cite[p.~84]{Newman}.
\end{proof}
\begin{rem}
The arguments of this section can be used to show that
if a function $f$ has $N$ tracts $D_1,\ldots,D_N$
and if $(s_k)_{k\geq 0}$ is a sequence 
in $\{1,\ldots,N\}$, then there exists $z\in D_{s_0}$ such 
that $f^k(z)\in D_{s_k}$ for all $k\in \N$ and
$f^k(z)\to\infty$ as $k\to\infty$.
Moreover, the set of all $z$ with this property has an unbounded
component.

If $f$ has infinitely many tracts $D_1,D_2,\ldots$, then
not every sequence $(s_k)$  in $\N$ is admissible.
This occurs, for example, with the function $f(z)=\exp(\exp z)$.
Noting that the tracts of $\exp(\exp z)$ are contained in
half-strips 
$\{z\in\C: \re z>0, (2j-1)\pi<\im z<(2j+1)\pi\}$
with $j\in\Z$,
we see that the question of when a sequence $(s_k)$ is admissible
is closely connected to the concept of an allowable 
itinerary considered in~\cite[Section~3]{Dev} and~\cite{DevKry}
and many other papers on exponential dynamics.
\end{rem}

\section{Meromorphic  functions with a direct tract}
\label{meromorphic}
Let $f$ be a transcendental function meromorphic in the
plane which has a direct tract $D$, and
let $\rho>\rho_0$ be as in the previous section.
We define
\begin{eqnarray*}
A(f,D)
&= &
\{z\in \C: \text{ there exists } L\in\N \text{ such that }
f^L(z)\in A(f,D,\rho)\}\\
&= &
\{z\in \C: \text{ there exists } L\in\N \text{ such that }
f^{n+L}(z)\in D  \\
& & \quad \text{ and }|f^{n+L}(z)| \geq M_D^n(\rho) \text{ for all } n\in\N \}.
\end{eqnarray*}
The set $A(f,D)$ is a variation of the set of fast escaping points 
of an entire function:
$$A(f)=\{z\in \C: \text{ there exists } L\in\N \text{ such that }
f^{n+L}(z)\geq M(R,f^n) \text{ for all } n\in\N \},$$
where $R$ is so large 
that $J(f)\cap D(0,R)\neq\emptyset$.
It is known~\cite{BH,RS05} that
$A(f)\neq\emptyset$, $\partial A(f)=J(f)$,
$A(f)\cap J(f)\neq\emptyset$ and $A(f)$ is completely invariant;
moreover all components of $A(f)$ are unbounded.
It can be shown that, in the definition of $A(f)$, we can replace
$M(R,f^n)$ by $M^n(R,f)$.
Note that for a given tract $D$, the set $A(f,D)$ may not be a subset
of $A(f)$; see Remark~2 at the end of this section.

Using arguments as in~\cite{BH,RS05} one can show that
the set $A(f,D)$ does not depend
on the choice of $\rho$, as long as $\rho>\rho_0$.
This justifies the notation where $\rho$ is suppressed.
It follows from Theorem~\ref{th4} and arguments given
in~\cite{RS05} that if $f$ is entire, then
all components of $A(f,D)$ are unbounded.

In this section and the next we 
prove a number of properties of $A(f,D)$ which are analogous to known 
properties of $A(f)$, when $f$ is entire, and also 
to similar properties of meromorphic functions with a finite number of 
poles; see~\cite{BH,Ere89,RS05,RS05b}.  
Many proofs in these two sections are modifications of proofs 
in~\cite{RS05b} and we do not always give full details.

We also define 
\[
Z(f) = \left\{z \in I(f): \tfrac{1}{n}\log\log |f^n(z)| \to \infty\mbox{ as } n \to \infty \right\}
\]
and we recall from~\cite[Lemma~7]{Ber93} that any periodic Fatou 
component of $f$ does not meet $Z(f)$. 

\begin{thm}\label{Aprops1}
Let $f$ be a {\tmf} with a direct tract~$D$.
Then the following properties hold:
\begin{itemize}
\item[(a)] $A(f,D)\ne \emptyset$ and, for each $z\in A(f,D)$, there exists
$L\in\N$ such that $f^L(z)$ lies in an
unbounded closed connected subset of $A(f,D)$;
\item[(b)] $A(f,D)$ is completely invariant under $f$;
\item[(c)] $A(f,D)\subset Z(f)$.
\end{itemize}
\end{thm}
\begin{proof}
Part~(a) follows from Theorems~\ref{th3} and~\ref{th4}, and the definition
of $A(f,D)$.

Part~(b) follows easily from the fact that $M_D(r)>r$ for $r\ge\rho$.

Part~(c) follows from part~(b) and~\eqref{growthB} which says
that $\log M_D(r)/\log r\to\infty$ as \mbox{$r\to\infty$.} 
\end{proof}
Our next result relates $A(f,D)$ to the Fatou set and Julia set of~$f$.
\begin{thm}\label{Aprops2}
Let $f$ be a {\tmf} with a direct tract~$D$.
Then the following properties hold:
\begin{itemize}
\item[(a)] if $U$ is a component of $F(f)$ such that 
$U\cap A(f,D)\ne \emptyset$, then $U\subset A(f,D)$
and $U$ is a wandering domain;
\item[(b)] $J(f)=\partial A(f,D)$;
\item[(c)] if $f$ has no wandering domains, then $J(f)=\overline{A(f,D)}$.
\end{itemize}
\end{thm}
\begin{proof}
To prove part~(a), let $U$ be a component of $F(f)$ which meets $A(f,D)$. 
By the complete invariance of $A(f,D)$ we can assume that $U$ meets 
$A(f,D,\rho)$, for some $\rho>0$, at a point $z_0$ say. Let $\Delta_0$ be 
a closed disc in $U$ with centre $z_0$. Then
\[|f^{n}(z_0)|>M_D^n(\rho)\quad\text{and}\quad f^n(z_0)\in D,\quad\text{for }n=0,1,\ldots.\] 
By~\cite[Theorem~3(a)]{RS00}, for example, there exists $C>0$ such that
for all $z\in \Delta_0,$ and $n=0,1,\ldots$, we have
\[
|f^n(z)|\ge |f^n(z_0)|^{1/C}>M_D^n(\rho)^{1/C}.
\]
Since $\log M_D^n(\rho)/\log M_D^{n-1}(\rho)\to\infty$ as $n\to \infty$, we 
deduce that, for all $z\in \Delta_0$ and all sufficiently large $n$,
\[
|f^n(z)|>M_D^{n-1}(\rho)>R\quad\text{so}\quad f^n(z)\in D,
\]
since $|f|=R$ on $\partial D$. Hence, for some $m\in\N$ we have 
$f^m(z)\in A(f,D,\rho)$. Thus $z\in A(f,D)$, 
so $\Delta_0\subset A(f,D)$. It follows by a compactness 
argument that $U\subset A(f,D)$, 
as required. Hence, by Theorem~\ref{Aprops1}(c), $U$ is a wandering domain.

The proof of part~(b) is identical to the proof of~\cite[Theorem~2(c)]{RS05b} 
and part~(c) follows immediately from parts~(a) and~(b).
\end{proof}

\begin{remark}
One consequence of Theorems~\ref{Aprops1} and~\ref{Aprops2} is that if $f$ 
is a transcendental meromorphic function with a direct tract $D$,
then none of the sets $J(f)$, $A(f,D)$, $Z(f)$ and $I(f)$ can contain
a free Jordan arc; see~\cite[Theorem~6]{RS05b} for details.
\end{remark}
\begin{remark}
Taking $R>e+1$ we see that $f(z)=\exp(-z)+\exp(\exp(z))$ has one
direct tract $D$ in the left half-plane and infinitely many 
direct tracts in the right half-plane. It is not difficult to
see that $A(f,D)\cap A(f)=\emptyset$. On the other hand,
denoting by $D'$ the direct tract that contains all large positive
real numbers, we see that $A(f,D')\subset A(f)$.

There are also examples of entire functions $f$ for which 
$A(f,D)\cap A(f)=\emptyset$ for every direct tract~$D$.
More specifically, there exists an entire function $f$ 
with exactly two tracts $D_1$ and $D_2$ such that 
$A(f,D_j)\cap A(f)=\emptyset$ for $j=1,2$. We only 
indicate very briefly how such a function can be constructed.
Let $\theta_j:(1,\infty)\to (0,\frac12 \pi)$ be continuous for $j=1,2$
and let $\Omega_1=\{re^{i\varphi}:r>1,|\varphi|<\theta_1(r)\}$
and $\Omega_2=\{re^{i\varphi}:r>1,|\varphi-\pi|<\theta_2(r)\}$.
There are several techniques for constructing
an entire function $f$  which has two tracts $D_1$ and $D_2$
which are ``close'' to $\Omega_1$ and $\Omega_2$,
provided the functions $\theta_1$ and $\theta_2$ are sufficiently ``nice''.
One such technique is the Kjellberg-Kennedy-Katifi
approximation method described in~\cite[Section~10.5]{Hay89},
another one is via Cauchy integrals 
(see~\cite[Section~VI.4]{Evg} or~\cite{Sta97}).
One also has control over the growth of~$f$ in the tracts
by choosing the sizes of
$\theta_j(r)$ as $r$ varies. In particular, one can arrange that
$M_{D_1}(r)$ is much bigger
than $M_{D_2}(r)$ in certain intervals and much smaller 
than $M_{D_2}(r)$ in other intervals. In this way it is possible 
to construct $f$ such that $M^n_{D_1}(\rho)$ and
$M^n_{D_2}(\rho)$ both grow much more slowly than $M^n(\rho,f)$ as 
$n\to\infty$. This then implies that $A(f,D_j)\cap A(f)=\emptyset$ 
for $j=1,2$.
\end{remark}

\section{Meromorphic functions and Baker wandering domains}
\label{functionwithbwd}
In this section we prove some results about {\bwd}s 
which, in particular, contain Theorems~\ref{continua} and
\ref{unboundedifnobwd} as special cases. Recall that if $U$ is a 
component of $F(f)$, then $U_n$ denotes the component of $F(f)$ 
such that $f^n(U)\subset U_n$, and that a Baker wandering domain 
is a component $U$ of $F(f)$ such that, for $n$ large 
enough, $U_n$ is a bounded multiply connected component of 
$F(f)$ which surrounds~0, and $U_n \to \infty$ as $n\to\infty$.
We use the notation $\widetilde{U}$ to denote the 
union of a set $U$ and its bounded complementary components. If 
$\widetilde{U}=U$, then we say that $U$ is {\it full}.  
\begin{thm}\label{bwd}
Let $f$ be a {\tmf} with a direct tract~$D$.
If $f$ has a {\bwd} $U$, then 
\begin{itemize}
\item[(a)] $\overline{U}\subset A(f,D)$, more precisely, there exist $N\in \N$ 
and $\rho>0$ such that, for $n\ge N$,
\[U_n\subset \wt{U_{n+1}}\quad\text{and}\quad \overline{U_n}\subset A(f,D,\rho);\]
\item[(b)] $\overline{U_n}\subset D$, for $n\ge N$, so $D$ is the only direct 
tract of $f$ and it has no unbounded complementary components;
\item[(c)] $A(f,D)\cap J(f)$ contains infinitely many bounded continua;
\item[(d)] $A(f,D)$ has exactly one unbounded component, as do $Z(f)$ 
and $I(f)$.
\end{itemize}
\end{thm}
Thus if $f$ has a direct tract $D$ and a Baker wandering domain, 
then $A(f,D)$ has a single unbounded component which must meet both
$F(f)$ and $J(f)$. We remark that there exist entire functions for which 
 $A(f,D)$ and $I(f)$ are connected, and $I(f)\subset J(f)$; see~\cite{Sta}.

Next we give a sufficient condition for {\bwd}s to exist.
\begin{thm}\label{bwdsufficient}
Let $f$ be a {\tmf} with a direct tract~$D$.
Then there is a constant $r_0>0$ such 
that if $U$ is a component of $F(f)$ which contains a Jordan curve 
surrounding  $\{z\in\C:|z|=r_0\}$, then $U$ is a {\bwd}.
\end{thm}
In Theorem~\ref{bwd}, we showed that 
$A(f,D)\cap J(f)$ contains continua
if $f$ has a {\bwd}. We now show that this conclusion also holds if 
$f$ does not have {\bwd}s.
\begin{thm}\label{nobwd}
Let $f$ be a {\tmf} with a direct tract~$D$
and suppose that $f$ has no {\bwd}s. Then
$A(f,D)\cap J(f)$ has at least one unbounded component, 
as do $Z(f)\cap J(f)$ and $I(f)\cap J(f)$.
\end{thm}
\begin{rem}
In particular, it follows from Theorem~\ref{nobwd}
that $J(f)$ has an unbounded component, and
it is this statement which will be proved first.
\end{rem}

To prove these results we use the concept of an outer sequence.
Starting from a 
full open set $G_0$ (not necessarily connected), we carry out the 
following process. First, form $B_0=G_0\cap D$, then take the image 
under $f$ and finally fill in any bounded complementary components to obtain 
the full open set
$G_1=\wt{f(B_0)}$. 
Then we repeat this process to obtain a sequence $(G_n)$ of full open sets, 
so long as this is possible. The following lemma shows that if $G_0$ is 
a large enough Jordan domain, then $(G_n)$ is
a sequence of Jordan domains whose 
boundaries tend to~$\infty$. If $G_0$ satisfies the hypotheses of 
Lemma~\ref{outerseq}, then we say that the sets $E_n=\C\setminus G_n$ 
form an {\it outer sequence} for $(f,D)$.
\begin{la}\label{outerseq}
Let $f$ be a {\tmf} with a direct tract~$D$ and let $M_D(r)$ be defined
by~\eqref{defM}.
Then there is a constant 
$r_0>0$ such that $D\cap \{z\in\C:|z|=r_0\}\ne\emptyset$ and if $G_0$ is a 
Jordan domain with boundary $\gamma_0$ that surrounds $\{z\in\C:|z|=r_0\}$, then 
the corresponding sets $G_n$ are Jordan domains with the following properties.
\begin{itemize}
\item[(a)] For $n=0,1,\ldots,$ the Jordan curve $\gamma_n=\partial G_n$ 
surrounds $\{z\in\C:|z|=r_n\}$, where
\[
 r_{n+1}>4M_D(r_n/4)>2r_n,\quad\text{for }n=0,1,\ldots.
\]
\item[(b)] For $n=0,1,\ldots,$
\[
\gamma_{n+1}\subset f(\gamma_n\cap D).
\]
\item[(c)] For $n=0,1,\ldots,$ any component of $f^{-1}(E_{n+1})$ which 
meets $E_n\cap D$ lies in $E_n\cap D$. 
\end{itemize}
\end{la}
\begin{proof}
Let $R$, $v$, $z_r$ and $F$ be as in Theorem~\ref{th2}, with $\beta=4$ and 
the corresponding value of $\alpha$. Then choose $r_0>0$ such that 
$D\cap \{z\in\C:|z|=r_0/4\}\ne\emptyset$,
\begin{equation} \label{5.1}
\text{each interval of the form } (r,2r), r\ge r_0/4,
\text{ contains a point outside } F,
\end{equation}
which is possible since $F$ has finite logarithmic measure, 
\begin{equation} \label{5.2}
\frac{\alpha}{a(r,v)}<1,\quad\text{for }r>r_0/4,
\end{equation}
which is possible since $a(r,v)\to\infty$ as $r\to\infty$,
\begin{equation} \label{5.3}
r_0>R,
\end{equation}
\begin{equation} \label{5.4}
 \mbox{the conclusion of Theorem~\ref{th2} holds for } r\in (r_0/4,\infty)\setminus F,
\end{equation}
\begin{equation} \label{5.5}
 M_D(r)>2r, \quad\mbox{for } r \ge r_0/4,
\end{equation}
which is possible since $B(r)/\log r=\log M_D(r)/\log r\to\infty$ as
$r\to\infty$.

Now by~\eqref{5.1} we can
choose $\rho_0\in (r_0/4,r_0/2)\setminus F$ and $\zeta_0$ such
that $\zeta_0=z_{\rho_0}$. Then $|f(\zeta_0)|=M_D(\rho_0)$ and, 
by~\eqref{5.4},
\[
f(A_0)\supset\left\{z\in\C:\tfrac14 M_D(\rho_0)<|z|<4M_D(\rho_0)\right\},
\]
where $A_0=D(\zeta_0,\alpha\rho_0/a(\rho_0,v))$. Recall from Theorem~\ref{th1}
that $A_0\subset D$. Also, $A_0\subset G_0$, since
$\alpha\rho_0/a(\rho_0,v)<\rho_0<r_0/2$, by~\eqref{5.2}.
Thus $A_0\subset B_0=G_0\cap D$,
so 
\[
f(B_0)\supset\left\{z\in\C:\tfrac14 M_D(\rho_0)<|z|<4M_D(\rho_0)\right\}.
\]
Clearly $f(B_0)$ is bounded and there are no components of $f(B_0)$ which
lie entirely outside the above annulus, since any such component would
have to meet the circle $\{z\in\C:|z|=R\}$. Thus $G_1=\wt{f(B_0)}$ is a Jordan
region and $\gamma_1=\partial G_1$ is a Jordan curve which lies outside
$\{z\in\C:|z|=r_1\}$, where $r_1=4M_D(\rho_0)$. Also,
by~\eqref{5.3} and~\eqref{5.5}, we have
\[r_{1}=4M_D(\rho_0)>4M_D(r_0/4)>2r_0>R.\]
Thus $\gamma_{1}\subset f(\gamma_0\cap D)$ because
$\partial f(B_0) \subset f(\partial B_0)$ and $|f|=R$ on
$\partial B_0\setminus \gamma_0$. 

Recall that, for $n=0,1,\ldots,$ the set $E_n$ denotes the complement of
$G_n$. If $K$ is any component of $f^{-1}(E_1)$ which meets
$E_0\cap D$, then $K$ must lie entirely in $D$ (because $|f|=R<r_1$
on $\partial D$) and cannot meet $B_0$ (since $f(B_0)\cap E_1=\emptyset$).
Thus $K$ is a subset of $E_0\cap D$.

We can now carry out this process with $\gamma_1$ and $r_1$ in place of
$\gamma_0$ and $r_0$, and continue repeatedly to produce the required
sequences $(G_n)$,
$(\gamma_n)$, $(r_n)$ and $(E_n)$ which satisfy parts~\mbox{(a)--(c)}. 
\end{proof}
\begin{rem}
The sequences $(G_n)$, $(\gamma_n)$ and $(E_n)$ defined in 
Lemma~\ref{outerseq} are
related to sequences of sets introduced in~\cite{RS05} and~\cite{RS05b}.
For a {\tef} $f$ the set
\[B(f)=\{z\in\C: \text{ there exists } L \in {\N} \text{ such that }
f^{n+L}(z)\notin \wt{f^n(\Delta)}, \text{ for } n \in {\N} \},\] 
where $\Delta$ is an open disc that meets $J(f)$, was defined 
in~\cite{RS05} and proved
equal to the set $A(f)$ defined earlier.
The above sequence $(G_n)$ is a
modification of the sequence $\left(\wt{f^n(\Delta)}\right)$. 
The definition of $B(f)$
was extended to {\tmf}s with a finite number of poles in~\cite{RS05b} using
the concept of an outer sequence.
\end{rem}

Outer sequences give us an alternative way to prove that points
are in $A(f,D)$. 
\begin{la}\label{A(f,D)criterion}
Suppose that $f$, $D$, $E_n$ and $r_n$ are as in Lemma~\ref{outerseq}.
If $z$ satisfies
\[f^n(z)\in E_n\cap D,\quad \text{for }n=0,1,2,\ldots,\]
then 
\[|f^n(z)|\ge 4M_D^n(r_0/4),\quad\text{for } n=0,1,2,\ldots,\]
so $z\in A(f,D,\rho)$, where $\rho=r_0/4$.
\end{la}
\begin{proof}
Since 
\[|f^n(z)|\ge r_n,\quad\text{for } n=0,1,2,\ldots,\]
and, by Lemma~\ref{outerseq}(a) and induction,
\[r_n\ge 4M_D^n(r_0/4),\quad\text{for } n=0,1,\ldots,\]
the result follows.
\end{proof}

\begin{rem}
We could use the outer sequences in Lemma~\ref{outerseq} to
define $B(f,D)$ in a similar way to the definition of $B(f)$ in~\cite{RS05b},
and then use Lemma~\ref{A(f,D)criterion} to prove that $B(f,D)=A(f,D)$.
\end{rem}

\begin{proof}[Proof of Theorem~\ref{bwd}]
Let $U_0=U$ be a Baker wandering domain. By the definition of a {\bwd}, we
deduce from Theorem~\ref{th4} and Theorem~\ref{Aprops2}(a) that
$U_n\subset A(f,D)$ for all sufficiently large $n$, so $U_0\subset A(f,D)$
by Theorem~\ref{Aprops1}(b). To prove that $\overline{U_0}\subset A(f,D)$,
we argue more carefully as follows. By the complete invariance of $A(f,D)$,
we can renumber $U_n$ so that
\[U_n\text{ surrounds }\{z\in\C:|z|=r_0\},\quad\text{for } n=0,1,\ldots,\]
where $r_0$ satisfies the hypotheses of Lemma~\ref{outerseq}. In particular,
for $n=0,1,\ldots$, we have $U_n\cap D\ne \emptyset$, by Lemma~\ref{outerseq},
and $|f(z)|>r_0>R$, for $z\in U_n$, so $\overline{U_n}\subset D$.  

Then take a Jordan curve $\gamma_0$ in $U_0$ surrounding 0 and define
$\gamma_n$, $G_n$ and $E_n$ as in Lemma~\ref{outerseq}
and $B_n=G_n\cap D$. Clearly
$\gamma_n\subset U_{n},$ for $n=0,1,\ldots$. We now show that 
\[U_n\subset\wt{U_{n+1}},\quad\text{for }n \ge 0.\] 
To prove this we show that, for $n\ge 0$, the curve $\gamma_{n+1}$ lies outside
$\gamma_n$. If this is false, then $\gamma_{n+1}$ lies inside $\gamma_n$ and
hence $G_{n+1}\subset G_n$. But then
\[G_{n+2}=\wt{f(B_{n+1})}\subset\wt{f(B_{n})}=G_{n+1}\subset G_n.\]
By induction, $G_m\subset G_n$, for $m>n$, which contradicts
Lemma~\ref{outerseq}(a). 

Thus $\overline{U_{n+1}}\subset E_{n}\cap D$, for $n=0,1,\ldots$,
so $\overline{U_{n+1}}\subset A(f,D,r_0/4)$, for $n=0,1,\ldots,$
by Lemma~\ref{A(f,D)criterion}, as required.

The fact that $\overline{U_n}\subset D$, for $n=0,1,\ldots$, implies
that $D$ is the only tract of~$f$ and all its complementary components are
bounded.

Part~(c) holds because $\partial U_n\subset A(f,D)\cap J(f)$,
for $n=0,1,\ldots$.

Part~(d) follows from part~(a), Theorem~\ref{th4} and Theorem~\ref{Aprops1}(c).
\end{proof}

\begin{proof}[Proof of Theorem~\ref{bwdsufficient}]
Let $r_0$ be the constant in Lemma~\ref{outerseq}. Denote the given
Jordan curve by $\gamma_0$ and define the corresponding outer sequence
$(E_n)$ as in Lemma~\ref{outerseq}. Then each $\gamma_n$ meets $A(f,D)$,
by Theorem~\ref{th4} and Theorem~\ref{Aprops1}(b). Thus the Fatou components $U_n$ are
wandering by Theorem~\ref{Aprops2}(a) and hence they are disjoint. Since
each component $U_n$ contains the curve $\gamma_n$ and does not meet any
$\gamma_m$ for $m\ne n$, it follows from Lemma~\ref{outerseq}(a)
that $U=U_0$ is a {\bwd}. 
\end{proof}

The proof of Theorem~\ref{nobwd} is longer. We begin by showing
that $J(f)$ has at least one unbounded component. First recall
from~\cite[Section~6]{RS05b} some general properties of the family $\J$
of components $J$ of $J(f)$. For $J\in \J$, we again use the notation $\wt{J}$
to denote the union of $J$ and its {\bcc}s, and we associate to each $J\in \J$
the set
$$
\Omega_J=\bigcup \left\{\wt{J'}:J'\in\J, \wt{J}\subset \wt{J'}\right\}.
$$
Note that $\partial\Omega_J$ is a subset of $J(f)$ since if
$z\in \partial \Omega_J$, then each open neighbourhood of $z$ must meet
some component of $J(f)$. Also, each $\Omega_J$ is
connected. The following lemma shows that each set $\Omega_J$ is full.

\begin{la}\label{closedset}
Let $f$ be meromorphic and for each $J\in\J$ let $\Omega_J$ be defined as above.
\begin{itemize}
\item[(a)] If $\Omega_J$ is bounded, then $\Omega_J$ is compact and full.
\item[(b)] If $\Omega_J$ is unbounded, then
\begin{itemize}
\item[(i)] either $\Omega_J=\wt{J'}$ for some unbounded $J'\in\J$
so $\Omega_J$ is closed and full,
\item[(ii)] or $\Omega_J$ is the union of a sequence of bounded sets of
the form $\wt{J_n}$, $n=0,1,\ldots,$ where $J_n\in \J$ and
$\wt{J_1}\subset \wt{J_2}\subset \cdots$, and
$\Omega_J$ is open and full.
\end{itemize}
\end{itemize}
\end{la}

Lemma~\ref{closedset} was proved in~\cite[Lemma~6 and its proof]{RS05b} under
the assumption that $f$ has a finite number of poles, but it was remarked
there that it is essentially a topological result about a closed set in $\C$.

Next we show that, with $f$ and $D$ as above, the sets $\Omega_J$ cannot
all be bounded.

\begin{la}\label{unboundedomega}
Let $f$ be a {\tmf} with a direct tract. Then $\Omega_J$ is unbounded
for at least one $J\in \J$.
\end{la}
\begin{proof}
Suppose that all $\Omega_J,\,J\in\J,$ are bounded.  Then all
$\Omega_J,\,J\in\J,$
are compact and full by Lemma~\ref{closedset}(a). Also note that
any two $\Omega_J$
are either disjoint or identical. We consider the set
\[
U=\C\setminus \bigcup \left\{\Omega_J:J\in\J\right\}.
\]
Then it is a topological result (see~\cite[proof of Lemma~7]{RS05b}) that
$U$ is nonempty, open and connected.

Now $U$ is evidently an unbounded subset of $F(f)$. Moreover $U$
must be a component of $F(f)$, indeed the only unbounded component
of $F(f)$ because any unbounded connected set must meet~$U$. If $D$ is a
tract of $f$ and $v$ is the corresponding subharmonic function defined
in~\eqref{definitionv}, then $v$ 
and hence $f$ has asymptotic value $\infty$ by Iversen's
theorem for subharmonic functions~\cite[Theorem~4.17]{HK1}. Therefore $U$ must
be invariant under $f$, since $f(U)$ cannot be
bounded. But we also know that $U\cap A(f,D)\ne \emptyset$, by
Theorem~\ref{th4}, so we obtain a contradiction to Theorem~\ref{Aprops2}(a).
\end{proof}

Next we relate Lemma~\ref{closedset} to the existence of {\bwd}s.

\begin{la}\label{bwdcriterion}
Let $f$ be a {\tmf} with a direct tract. Then $f$ has a {\bwd} if and
only if Lemma~\ref{closedset}(b)(ii)
holds with $\Omega_J=\C$, for some $J\in\J$.
\end{la}
\begin{proof}
It is clear that if $f$ has a {\bwd}, then $\C$ is 
the union of bounded sets of the form $\wt{J_n}$, $n=0,1,\ldots,$ 
where 
\begin{equation} \label{Jcond}
J_n\in \J \quad\text{and}\quad \wt{J_1}\subset \wt{J_2}\subset \cdots.
\end{equation} 

Suppose, on the other hand, that, for some $J\in\J$, the set
$\Omega_J=\C$ is the union of bounded sets of the form $\wt{J_n},\, n=1,2,\ldots,$
such that~(\ref{Jcond}) holds.

Choose $N$ so large that $\{z\in\C:|z|=r_0\}\subset \wt{J_N}$, where $r_0$ satisfies the
hypotheses of Lemma~\ref{outerseq}. Then $E=J(f)\cup \{\infty\}$ is
a closed subset of $\hat{\C}$ having $J_N$ and $J_{N+1}$ amongst its
components. Thus there is a simple polygon $\gamma_0$ separating $J_N$ and $J_{N+1}$,
and lying in $F(f)$. Since $\gamma_0$ surrounds $\{z\in\C:|z|=r_0\}$, we
deduce from Theorem~\ref{bwdsufficient} that the component $U$ of
$F(f)$ which contains $\gamma_0$ is a {\bwd}, as required.
\end{proof}

It follows from Lemmas~\ref{closedset},~\ref{unboundedomega}
and~\ref{bwdcriterion} that if $f$ and $D$ are as above and
$f$ has no {\bwd}s, then for some $J\in\J$ the set $\Omega_J$
is unbounded and either $\Omega_J=\wt{J'}$, for some $J'\in\J$, or
$\Omega_J$ is a union of bounded sets of the form $\wt{J_n}, \,
J_n\in \J,$ and $\Omega_J \ne \C$. In the former case, $J(f)$
certainly  has an unbounded component, namely $J'$. In the latter
case, $\Omega_J$ is unbounded, open, connected and full. As
in~\cite[page~240]{RS05b}, we then deduce that the boundary of any
complementary component of $\Omega_J$ is an unbounded subset of $J(f)$,
as required. 
Thus in either case $J(f)$ has an unbounded component.

To complete the proof of
Theorem~\ref{nobwd}, we need another topological lemma.

\begin{la}\label{preimage}
Let $f$ be a transcendental meromorphic function
with a direct tract $D$ and suppose that $(E_n)_{n\geq 0}$ is
an outer sequence for $(f,D)$, given by Lemma~\ref{outerseq}, with
corresponding Jordan curves $\gamma_n$. For $n\ge 1,$ let $K_n$ be
a closed subset of $J(f)\cap E_n$ which meets $\gamma_n$ and
whose components are unbounded. Then
$K_{n-1}=f^{-1}(K_n)\cap E_{n-1}\cap D$ is a closed subset of
$J(f)\cap E_{n-1}\cap D$
which meets $\gamma_{n-1}\cap D$ and whose components are unbounded.
\end{la}

\begin{proof}
First, $K_{n-1}=f^{-1}(K_n)\cap E_{n-1}\cap D$ is a
closed subset of $J(f)\cap E_{n-1}\cap D$, by the complete invariance of
$J(f)$, the continuity of $f$, and the fact that
$f^{-1}(E_n)\cap\,\partial D=\emptyset$. Also
$K_{n-1}$ meets $\gamma_{n-1}\cap D$, by Lemma~\ref{outerseq}(b). Now let
$C$ be a component of $K_{n-1}$. Then $f(C)$ is connected, so it is
a subset of a component $\Gamma$ of $K_n$.
The component of $f^{-1}(\Gamma)$ which contains $C$ is itself
contained in $K_{n-1}$, by Lemma~\ref{outerseq}(c), so we deduce that $C$
is a component of $f^{-1}(\Gamma)$ which lies in
$E_{n-1}\cap D$. Because $\Gamma$ is unbounded by hypothesis
and $f$ is holomorphic in~$D$, the set $C$ is unbounded.
This completes the proof of Lemma~\ref{preimage}.
\end{proof}

To prove that $A(f,D)\cap J(f)$ has at least one unbounded
component, we apply Lemma~\ref{preimage} repeatedly. First we may assume that
the outer sequence $(E_n)$ is chosen in such a way that each
$\gamma_n,\,n=0,1,\ldots,$ meets an unbounded component of $J(f)$.
Then, for $n=0,1,\ldots,$ let $K^0_n$ denote the union of {\it
all} the unbounded components of $J(f)\cap E_n$. The set $K^0_n$
is nonempty and meets $\gamma_n$, and we can show that $K^0_n$ is a
closed set as in~\cite[proof of Lemma~9]{RS05b}.

We deduce, by Lemma~\ref{preimage}, that 
$K^1_{n-1}=f^{-1}(K^0_n)\cap E_{n-1}\cap D$ is a closed subset of
$J(f)\cap E_{n-1}\cap D$, which meets $\gamma_{n-1}$ and whose components
are unbounded. Applying Lemma~\ref{preimage} repeatedly in this way 
we obtain, for $k=0,1,\ldots,n,$ a closed
subset $K^k_{n-k}$ of $J(f)\cap E_{n-k}\cap D$, which meets $\gamma_{n-k}$
and whose components are unbounded. Clearly
\[
J(f)\cap E_0\cap D\supset K^0_0\supset K^1_0\supset K^2_0 \supset \cdots.
\]
Then $\left(K^n_0\cup \{\infty\}\right)_{n\geq 0}$
is a nested sequence of continua in $\hat{\C}$,
each meeting $\gamma_0$ and including $\infty$, whose intersection, $K$ say, must be a continuum in
$\hat{\C}$ meeting $\gamma_0$ and including $\infty$. Thus any component $\Gamma$ of
$\bigcap _{n=0}^{\infty} K^n_0=K\setminus \{\infty\}$ is unbounded.

But if $z\in \Gamma$, then we have $z\in K^{n+1}_0$ so $f^n(z)\in
J(f)\cap E_n\cap D$, for $n=0,1,\ldots.$ Hence $z\in A(f,D)\cap J(f)$,
by Lemma~\ref{A(f,D)criterion}. Thus $\Gamma\subset A(f,D)\cap J(f)$, so
the proof of Theorem~\ref{nobwd} is complete.

\section{Functions with a logarithmic tract}
\label{classB}

Let $f$ be a
function with a direct tract~$D$. We consider the set
\begin{eqnarray*}
A'(f,D) &= & \{z\in D: \text{ there exists } L\in\N \text{ such
that }
f^{n}(z)\in D \\
& & \quad  \text{ and } |f^{n+L}(z)| \geq M_D^n(\rho) \text{ for
all } n\in\N \},
\end{eqnarray*}
whose definition does not require that $f$ is defined outside 
of~$D$. As with the set $A(f,D)$ that we considered earlier, the set
$A'(f,D)$ does not depend on $\rho$, as long as $\rho$ is chosen
to be sufficiently large. The reason that we considered the set
$A(f,D)$ earlier (rather than the set $A'(f,D)$) was that, for an
entire function, all the components of $A(f,D)$ are unbounded
whereas the set $A'(f,D)$ may have bounded components;
we give an example of such a function at the end of this
section. In the current situation, however, the function $f$ may
not be defined outside of $D$ and so it makes sense to consider
the set $A'(f,D)$. Note that, for a meromorphic function $f$ with 
a direct tract $D$ and for $\rho$ large enough, we have 
\begin{equation} \label{AinI}
A(f,D,\rho)\subset A'(f,D) \subset A(f,D)\subset I(f).
\end{equation}

The main results of this section give estimates on the size of the 
set $A'(f,D)$ when $D$ is a logarithmic tract of~$f$.

\begin{thm}\label{box2}
Let $D$ be a logarithmic tract of~$f$. Then the upper box
dimension of $A'(f,D)$ is equal to~$2$.
\end{thm}

\begin{thm}\label{Hausdorff}
Let $D$ be a logarithmic tract of~$f$. Then the Hausdorff
dimension of $\overline{A'(f,D)}$ is strictly greater than~$1$.
\end{thm}

We also note that the following result can be deduced from
Lemma~\ref{nderiv} below as in the proof of~\cite[Theorem~A]{RS99a}.

\begin{thm}\label{J(f)}
Let $f$ be a transcendental meromorphic function with a
logarithmic tract~$D$. Then $A(f,D) \subset J(f)$.
\end{thm}

It follows from Theorem~\ref{box2}, Theorem~\ref{J(f)} and~(\ref{AinI}) 
that if $f$ is a transcendental meromorphic function with a logarithmic
tract then the upper box dimension of $I(f)\cap J(f)$ (and so that of 
$J(f)$) is equal to~2. Hence, by~\cite[Theorem~1.2]{RS05c}, 
the packing dimension of $J(f)$ is also equal to~2. This
proves Theorem~\ref{box}. On the other hand, Theorem~\ref{hausdorff}
follows from Theorem~\ref{Hausdorff}, Theorem~\ref{J(f)} and
the fact that $J(f)$ is closed.

While our proof of Theorem~\ref{hausdorff}
is based on the set of escaping points, Baranski, Karpinska and
Zdunik~\cite{BKZ} have recently shown that the set of points in the
Julia set whose orbit is bounded and contained in a logarithmic
tract also has Hausdorff dimension greater than 1, thereby giving
an alternative proof of Theorem~\ref{hausdorff}.

We now state a number of properties of functions with a
logarithmic tract. Most of these results can be proved in a
similar way to the analogous results for transcendental
meromorphic functions in the class~$B$. We then explain how
Theorem~\ref{box2} and Theorem~\ref{Hausdorff} can be proved by
making small changes to the proofs for meromorphic
functions with finitely many poles in the class~$B$.

First we note the following result.
\begin{la}\label{bdy}
Let $D$ be a logarithmic tract of~$f$.
Then $D$ is simply connected.
If  $|f(z)|=R$ for $z\in\partial D$ and $r>R$,
then $\{z\in D: |f(z)|>r\}$ is a logarithmic tract of~$f$, which
is bounded by an unbounded simple analytic curve.
\end{la}
The following result can be deduced from Koebe's distortion
theorem as in~\cite[Lemma~2.3]{S96} by taking a suitable cover of
the annulus by discs. Here $L = 81$.
\begin{la}\label{annulus}
Let $D$ be a logarithmic tract of $f$ with $|f(z)| = R$ for $z \in
\partial D$. Let $|w_0| = r \geq 2R$ and let $g$ be a branch of
$f^{-1}$ defined at $w_0$ such that $g(w_0) \in D$. If $g$ is
analytically continued along a curve $\gamma$ that winds at most
once round 0 and lies within $\{w\in\C: 4r/5 \leq |w| \leq 5r/4\}$,
then
\[
 |g'(w_0)|/L^{26} \leq |g'(w)| \leq L^{26}|g'(w_0)|,
\]
for each $w \in \gamma$.
\end{la}
The next result can be proved in a similar way to the analogous
result for meromorphic functions in the class $B$
by using Lemma~\ref{bdy}; see~\cite[Lemma~2.2]{RS99}, for example.
\begin{la}\label{deriv}
Let $D$ be a logarithmic tract of~$f$. There exists $R_1(f)
>0$ such that if $z \in D$ and $|z|,|f(z)|>R_1(f)$, then
\[
  |f'(z)| > \frac{|f(z)|\log |f(z)|}{16 \pi |z|}.
\]
\end{la}
By repeatedly applying Lemma~\ref{deriv}, we obtain the following.
\begin{la}\label{nderiv}
Let $D$ be a logarithmic tract of $f$ and let $n\in\N$. If $f^k(z) \in D$ and
$|f^k(z)|>\max \{ R_1(f), e^{16 \pi}\}$, for $0 \leq k \leq n$,
then
\[
  |(f^n)'(z)| > \frac{|f^n(z)|\log |f^n(z)|}{16 \pi |z|}.
\]
\end{la}
Using Koebe's distortion theorem together with Lemma~\ref{deriv},
the following result can be proved in the same way as the
analogous result for entire functions in the class~$B$;
see~\cite[Lemma~2.6]{S96}.
\begin{la}\label{iterates}
Let $D$ be a logarithmic tract of~$f$. There exists $R_2(f) \geq
R_1(f)$ such that, if $f^k(z) \in D$ for $0 \leq k <n$ and
$|f^k(z)| \geq R_2(f)$ for $0 \leq k \leq n$, then the branch $g$
of $f^{-n}$ that maps $f^n(z)$ to $z$ is defined on
$D(f^n(z),|f^n(z)|/4)$ and satisfies, for each $0 \leq k <n$:
\begin{enumerate}
\item for each $K \geq 4$,
\[
  f^kg(D(f^n(z),|f^n(z)|/K)) \subset D(f^k(z),|f^k(z)|/(4K));
\]

\item $f^k \circ g$ is univalent in $D(f^n(z),|f^n(z)|/4)$;

\item if $w \in D(f^n(z),|f^n(z)|/8)$ then
\[
  |(f^k \circ g)'(f^n(z))|/L \leq |(f^k \circ g)'(w)| \leq L|(f^k \circ g)'(f^n(z))|.
\]
\end{enumerate}
\end{la}

The next lemma gives information on the growth in a
logarithmic tract. It can be proved in the same way as the
analogous result for entire functions in the class~$B$;
see~\cite[Lemma~3.5]{RS05c} and also the similar argument we 
give at the end of Section~\ref{proofgrowth}.

\begin{la}\label{growthla}
Let $D$ be a logarithmic tract of~$f$. There exist constants
$c=c(f)
>0$ and $R_3 = R_3(f)>1$ such that
\[
\log M_D(r) \geq cr^{1/2} \geq (\log 2r)^2, \quad\text{for } r \geq R_3.
\]
\end{la}

We use the following result in the proof of
Theorem~\ref{Hausdorff};
this requires a slightly different proof to that given for the
analogous result in~\cite{RS06}. We use $C(r)$ to denote
the circle
$\{z\in\C:|z|=r\}$.

\begin{la}\label{curve}
Let $D$ be a logarithmic tract of~$f$. There exists $R_4(f)>0$
such that, if $r\geq R_4(f)$, then there is an unbounded simple
analytic curve $\Gamma \subset D$ with $|f(z)|=r$ on $\Gamma$ and
$\Gamma \cap C(r) \neq \emptyset$.
\end{la}

\begin{proof} If $r$ is sufficiently large, then 
$M_D(r)>r$.
For such $r$, we take a point $z \in D \cap C(r)$ for which
$|f(z)| = M_D(r)$ and let $D_r$ denote the component of
$f^{-1}(\C\setminus\overline{D}(0,r))$
which contains~$z$. 
Lemma~\ref{bdy} implies that $D_r$ is a logarithmic tract
whose  boundary $\Gamma$ is an unbounded simple analytic curve on which
$|f(z)|=r$. 
We have seen that $D_r \cap C(r) \neq \emptyset$ and,
since $D_r$ is simply connected, it follows that $\Gamma \cap
C(r) \neq \emptyset$.
\end{proof}

The following result will also be used in the proof of
Theorem~\ref{Hausdorff}. This can be proved in the same way as
the analogous result for meromorphic functions with finitely many
poles in class $B$;
see~\cite[Lemma~2.9]{RS06}. (Note that by using Lemma~\ref{curve}
we are able to ensure that the pre-images of the annulus lie in
the tract~$D$.)

\begin{la}\label{pre-images}
Let $D$ be a logarithmic tract of~$f$. There exists $R_5(f)>0$
such that, if $|w_0| = r \geq R_5(f)$, then there exists $z_0$
such that
 \begin{enumerate}
\item $f(z_0) = w_0$;

\item $299r/300 \leq |z_0| \leq 301r/300$;

\item $D(z_0,r/100) \cap D$ contains at least two pre-images under 
$f$ of each point in $\{w\in\C: 3r/4 < |w| < 5r/4\}$.
 \end{enumerate}
\end{la}

In order to prove Theorem~\ref{box2}, we also need an estimate on
the growth of the function $a(r,v)$ which was defined in~(\ref{defarv}).

To do this we use the following well-known lemma about real
functions, whose short proof we include for completeness.
\begin{la} \label{growthlemma4}
Let $x_0,\eps>0$ and let $h,H: [x_0,\infty)\to (0,\infty)$, with
$h$ increasing and $H(x)=\int_{x_0}^x h(s)ds+H(x_0)$. Then
there exists a set $E\subset [x_0,\infty)$ of finite measure such
that
\begin{equation} \label{hH} 
h(x)\leq H(x)^{1+\eps},
\end{equation}
for $x\geq x_0$, $x\notin E$.
\end{la}
\begin{proof}
Defining $E$ as the set of all $x\geq x_0$ where~\eqref{hH} does
not hold, we have
$$\int_{E} dx\leq \int_{E} \frac{h(x)}{H(x)^{1+\eps}}\,dx
\leq \int_{x_0}^\infty \frac{h(x)}{H(x)^{1+\eps}}\,dx =
\int_{H(x_0)}^\infty \frac{1}{u^{1+\eps}}\,du<\infty$$ and the
conclusion follows.
\end{proof}
For a subharmonic function $v$ we can take
$h(x)=a(e^x,v)$ and $H(x)=B(e^x,v)$ in Lemma~\ref{growthlemma4}.
With $r=e^x$ this yields the following result.
\begin{la} \label{growthlemma5}
Let $v:\C\to [-\infty,\infty)$ be subharmonic and let $\eps>0$.
Then there exists a set $F\subset [1,\infty)$ of finite
logarithmic measure such that
\begin{equation} \label{aB} 
a(r,v)\leq B(r,v)^{1+\eps},
\end{equation}
for $r\geq 1$, $r\notin F$.
\end{la}
We are now in a position to show how Theorem~\ref{box2} can be
proved by modifying the proof of~\cite[Theorem~1.1]{RS05c}. The
following changes are needed for this proof to be applied to a
function $f$ with a logarithmic tract $D$ instead of to an entire
function in the class~$B$.

1. In~\cite{RS05c} we showed that the upper box dimension of the
set $A(f)$ is equal to 2 when $f$ is an entire function in the
class~$B$. The set $A(f)$ should be replaced by the set $A'(f,D)$
throughout. The proof of~\cite[Theorem~1.1]{RS05c} uses the fact
that the set $A(f)$ contains an unbounded closed connected set
--- Theorem~\ref{th4} and~(\ref{AinI}) show that the set $A'(f,D)$ 
also has this property.

2. In~\cite[Section 3]{RS05c} we gave several properties of
entire functions in the class~$B$. These properties also hold
for functions with a logarithmic tract; see Lemmas~\ref{deriv},~\ref{nderiv},~\ref{iterates} and~\ref{growthla}.

3. The argument in~\cite[Section 5]{RS05c} develops the method of 
Eremenko based on Wiman-Valiron theory in order to construct a
fast escaping orbit $(f^n(z_0))$ such that $|f^{n+1}(z_0)| \sim
M(|f^n(z_0)|,f)$. The Wiman-Valiron estimates (5.3)-(5.5)
in~\cite{RS05c} can be replaced by~\eqref{wvlike} 
and~\eqref{wvlike3}. 
Throughout the argument, the central
index $N(r)$ should be replaced by $a(r,v)$ and the maximum
modulus $M(r,f)$ should be replaced by $M_D(r)$ to give an orbit
$f^n(z_0)$ lying in $D$ such that $|f^{n+1}(z_0)| \sim
M_D(|f^n(z_0)|)$. In ~\cite[Section 5]{RS05c}, Lemma~5.1 is used
to establish (5.10) and (5.11) which are then used to establish
(5.16). It follows from~\eqref{comparisonBa} and Lemma~\ref{growthlemma5} that
the sequence $(r_n)$ can be chosen so that the inequalities of
(5.16) hold with $N(r)$ replaced by $a(r,v)$. (As
in~\cite{RS05c}, there is plenty of freedom in choosing the
sequence $(r_n)$ and so a set of finite logarithmic measure can be
avoided.)

3. In the current context, the inverse function $G$ which appears
in the proof in~\cite{RS05c} is defined in a similar way in
relation to points on the orbit $\left(f^n(z_0)\right)$, but now its image
values lie entirely in the tract.

We now show how Theorem~\ref{Hausdorff} can be
proved by modifying the proof of~\cite[Theorem~3]{RS06}. The
following changes are needed for this proof to be applied to a
function $f$ with a logarithmic tract $D$ instead of to a
meromorphic function with finitely many poles in the class~$B$.

1. In~\cite[Section 2]{RS06} we gave several properties of
functions with finitely many poles in the class~$B$. The analogous
results for functions with a logarithmic tract were stated at the
beginning of this section.

2. In~\cite[Section 3]{RS06} we studied a set $I$ consisting of
certain pre-images of a point $z_0$ in $J(f)$. Here we insist (as
we may) that all these pre-images lie in the logarithmic tract~$D$.
Also, we choose the point $z_0$ to lie in $A'(f,D)$.

3. In~\cite[Lemma 3.1]{RS06}, for each $n \in \N$, we constructed
$2^n$ simple curves $\gamma_{n,i}$. By Lemma~\ref{curve} and
Lemma~\ref{pre-images}, we can choose each such curve so that it
begins at a point in~$D$. Since the image of the curve lies
outside a large disc, it follows that the whole curve lies in~$D$.
As there are no poles in $D$ and the image of the curve is
unbounded, the curve must join a point in $D$ to infinity as
required.

\begin{rem}
As mentioned at the beginning of this section, there exists 
an entire function with a direct tract $D$ such that 
$A'(f,D)$ has a  bounded component.
An example with this property is given by
$$f(z)=e^{z+2}-e-1$$
and $D=\{z\in \C:|f(z)|>R\}$ if $R>0$ is sufficiently small. We
note that $D$ is connected and hence a tract. It follows that the
sets $A(f,D)$ and $A(f)$, which were introduced in Section 3, are
equal and have no bounded components. To see that $A'(f,D)$ has
bounded components we note that $f(-1)=-1$ while $f(x)>x$ for
$x>-1$. It is not difficult to deduce from this that
$(-1,\infty)\subset A(f,D)$. In fact, noting that $f$ is
conjugate to the function $z\mapsto \lambda e^z$ for
$\lambda=e^{1-e}$ it can be deduced from the results
in~\cite{Dev,DevKry} that $(-1,\infty)$ is a component of
$A(f,D)$.

Now $f$ has a zero $\xi\in (-1,0)$ and  there exists a unique
$\eta\in (\xi,\infty)$ with $f(\eta)=R$. It follows that
$(\eta,\infty)\subset A'(f,D)$, and since $A'(f,D)\subset A(f,D)$
this implies that  $(\eta,\infty)$ is a component of $A'(f,D)$.
However, if $R$ is small enough then there exists $x\in D\cap
(-1,\xi)$ satisfying $f(x)\in (\eta,\infty)$ and $|f(x)|>R$. It follows that
$x\in A'(f,D)$, but the component of $A'(f,D)$ which contains $x$
is contained in $(-1,\xi)$ and hence is bounded.
\end{rem}

\section{Results of Teichm\"uller and Selberg}
\label{teichm}
It seems reasonable to expect that the dynamics of
a meromorphic function with ``few poles'' are 
similar to that of an entire function. (In fact it was this 
question that motivated our work.) We will state two results of
Teich\-m\"uller and Selberg which imply that
this is indeed the case 
under suitable additional hypotheses.

Teich\-m\"uller's result~\cite{Tei} says that a meromorphic function
in class $B$ has a logarithmic singularity if it has few poles and
if the multiplicity of the poles is uniformly bounded. Using the
standard terminology of Nevanlinna theory~\cite{Hay64,Nev53},
Teich\-m\"uller's result
can be stated as follows.
\begin{prop} \label{teichmueller}
Let $f$ be a transcendental meromorphic function in the class $B$, and
suppose that there exists $N\in \N$ such that the poles of $f$
have multiplicity at most~$N$. 
If $\delta(\infty,f)>0$ or, more generally, if $m(r,f)$ is unbounded,
then $f$ has a logarithmic singularity over infinity.
\end{prop}
Teichm\"uller states this only for the case where the singularities 
of the inverse function of $f$ lie over finitely many points.
(This class of functions is usually denoted by $S$ today.)
However, his proof extends to class $B$ without change.

Selberg~\cite{Sel45} gave a variation of Teichm\"uller's result
which even includes some functions that are not in the class~$B$.
His proof yields the following result.
\begin{prop} \label{selberg}
Let $f$ be a transcendental meromorphic function.  Suppose that 
there exist $R>0$ and $N\in\N$ such that for each component U of
$f^{-1}(\CC\setminus \overline{D}(0,R))$
which contains a pole the map
$f:U \to \CC\setminus \overline{D}(0,R))$
is a proper map of degree at most~$N$.
If $\delta(\infty,f)>0$, or more generally, if $m(r,f)/\log r$ is unbounded,
then $f$ has a direct singularity over infinity.
\end{prop}
Note that functions in class $B$ whose poles have multiplicity
at most $N$
satisfy the hypothesis of Selberg's result.
So a slightly weakened form of Teichm\"uller's result follows from
Selberg's result: the hypothesis that $m(r,f)$ is unbounded has to be
replaced by the hypothesis that $m(r,f)/\log r$ is unbounded.

\section{Examples}
First we remark that it is easy to find a transcendental entire
function which has a logarithmic tract but is not in class~$B$.
For example, the functions $f(z)=\exp(z^2)\cos(\sqrt{z})$
and $f(z)=\exp(z^2)/\Gamma(-z)$ 
are not in class~$B$ but each has a logarithmic tract which 
includes most of the negative real axis and a direct (non-logarithmic)
tract which includes most of the positive real axis (apart from
gaps near the the zeros of $f$). 

Here we give several examples, each with a direct tract and 
infinitely many poles.
\label{examples}
\begin{ex} \label{ex1}
Let $\lambda\in\C\setminus\{0\}$ and 
$$f(z)=\lambda\frac{e^{2z}-1}{e^z-1/z}.$$
Then $|f(z)|\leq |\lambda|(e^2+1)/(e-1)$ for $\re z=1$ while 
$|f(x)|\to\infty$ as $x\to\pm\infty$ for $x\in\R$.
For $R>|\lambda|(e^2+1)/(e-1)$ we thus find that the set
$f^{-1}(\CC\setminus \overline{D}(0,R))$ has at least two
unbounded components, one containing the large positive numbers and 
one containing the large negative numbers. The component containing
the large positive numbers is a direct tract, while the 
one containing the large negative numbers is not.
Now, $f$ is bounded on each line $\{z\in\C:\re z=a\}$ for $a\geq 1$ 
and it is not difficult to check that $f'$ has no zeros in 
$\{z\in\C:\re z\geq a\}$ for large enough $a>0$. Also $f$ has 
no finite asymptotic values. Thus, for large enough $R>0$ the direct tract containing
the large positive numbers is actually a logarithmic tract.

For all $\lambda\in\C\setminus\{0\}$, the function $f$ 
has a superattracting fixed point at~$0$. For $\lambda>1$
the function has a Baker domain containing all large negative numbers.
For $\lambda=1$ there are infinitely many Baker domains in the left
half-plane; see~\cite[Theorem~1]{RS99}.
Figure~\ref{fig1} shows the Julia and Fatou sets for the parameters
$\lambda=\frac12$, $\lambda=1$ and $\lambda=2$
(from left to right). The attracting basin
of zero is drawn grey, the points tending to infinity in the 
logarithmic tract are black and the remaining points (including
those in the Baker domains) are white.
The range shown is given by $-10\leq \re z\leq 8$, $|\im z|\leq 12$.
\begin{figure}[htb]
\begin{center}
\includegraphics[width=4.5cm]{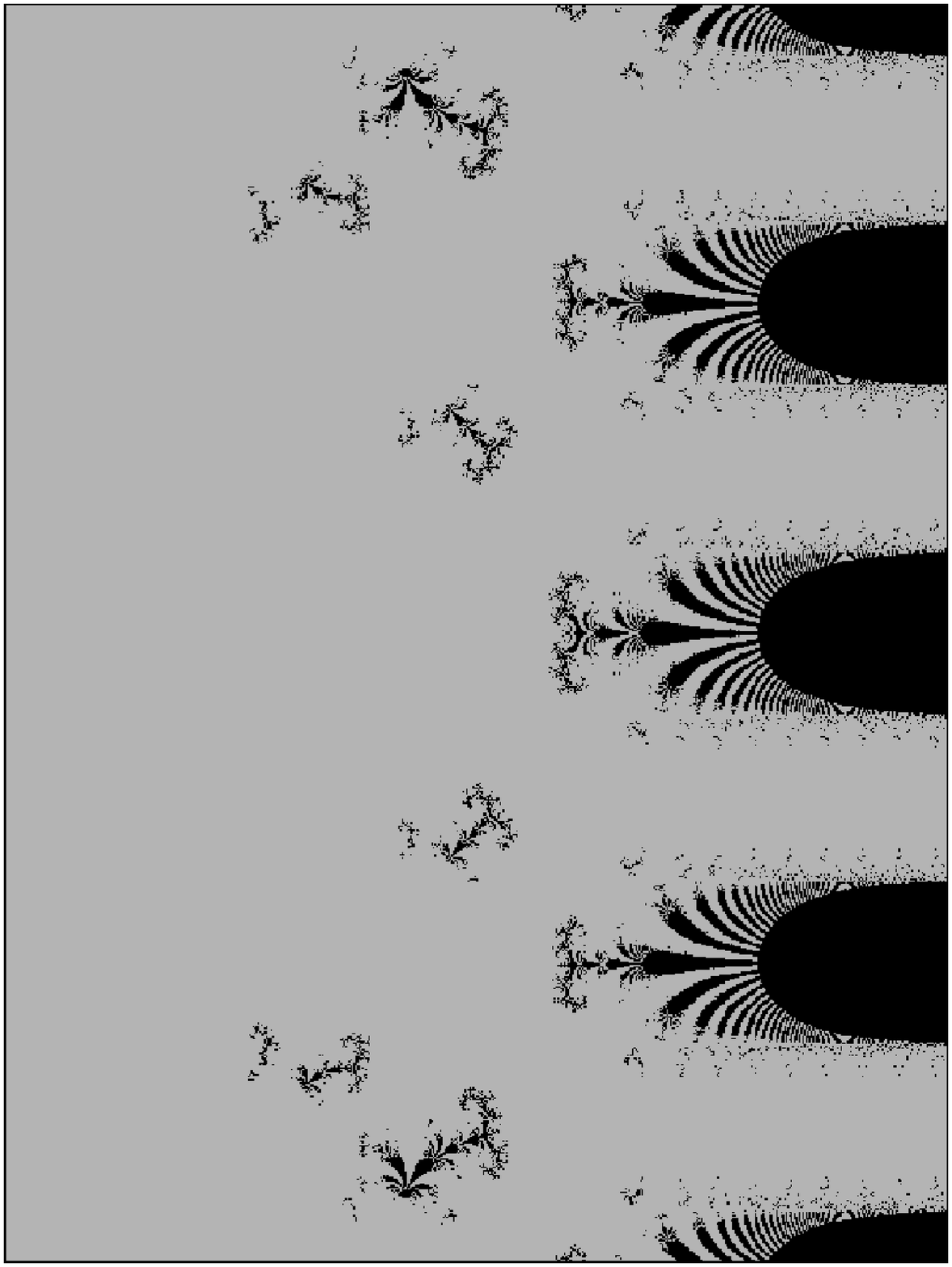}
\quad
\includegraphics[width=4.5cm]{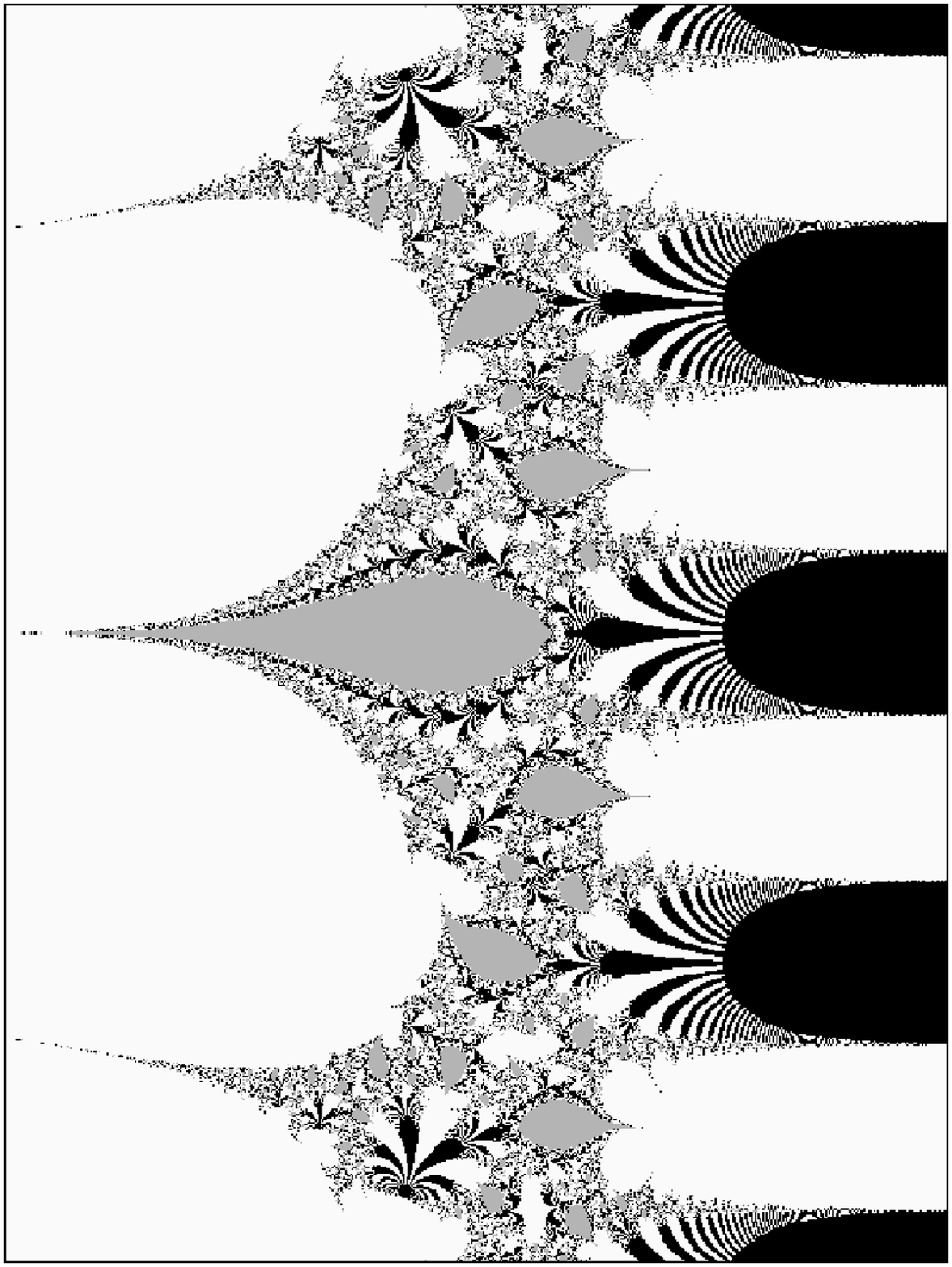}
\quad
\includegraphics[width=4.5cm]{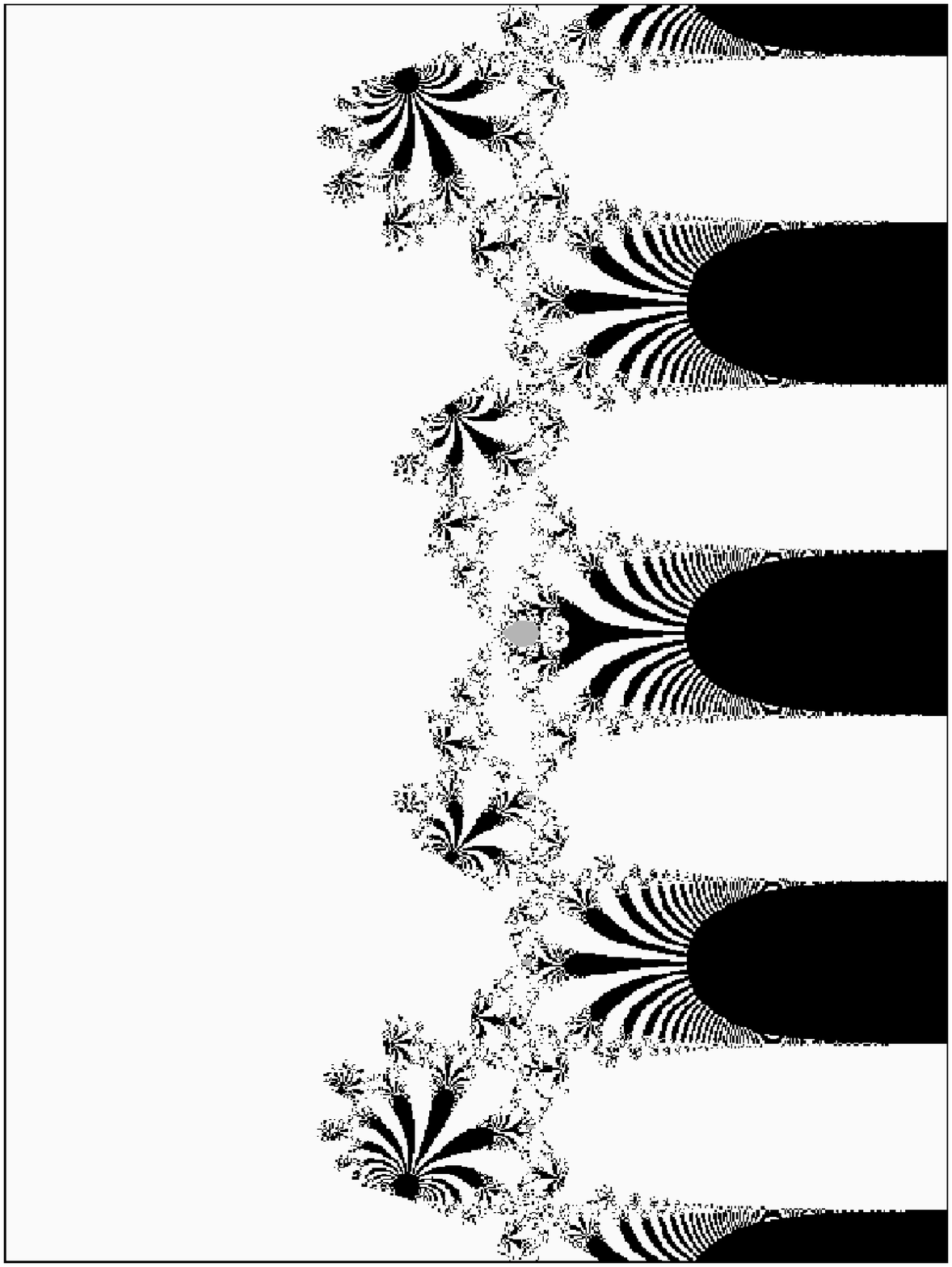}
\caption{The Julia sets of the functions from Example~\ref{ex1}.}
\label{fig1}
\end{center}
\end{figure}
\end{ex}
\begin{ex} \label{ex2}
It follows from Stirling's formula that the gamma function 
$\Gamma$ has a direct tract in the right half-plane.
Of course, this also holds for $\Gamma(z+a)$ with $a\in\C$.
One can also show that $\Gamma(z)$ and hence the functions $\Gamma(z+a)$
are in class $B$ and thus the direct tract is actually logarithmic.
We note that the gamma function also illustrates the results
of Section~\ref{teichm}, since $\delta(\infty,\Gamma)=1$.

The function $\Gamma(z+1)$ has an attracting fixed point at~$1$.
The point~$0$ is in the basin of attraction of this fixed point,
and since the gamma function is close to $0$ in the left half-plane
except for small neighbourhoods of the poles, large parts of the 
left half-plane are also contained in the attracting basin of~$1$.

The function $\Gamma(z)$ also has an attracting fixed point at~$1$.
However, $0$ is not in the basin of attraction, but $0$ is a pole.
This explains the rather complicated structure of the Julia set
in the left half-plane, which is mapped into a neighbourhood of~$0$.

Figure~\ref{fig2} shows the Julia and Fatou sets of
$\Gamma(z)$ (left) and $\Gamma(z+1)$ (middle).
The attracting basin of the fixed point at~$1$ is 
drawn white, the points tending to infinity in the 
logarithmic tract are black and the remaining points
are grey.
The range shown is given by $-5\leq \re z\leq 10$, $|\im z|\leq 10$.
\begin{figure}[htb]
\begin{center}
\includegraphics[width=4.5cm]{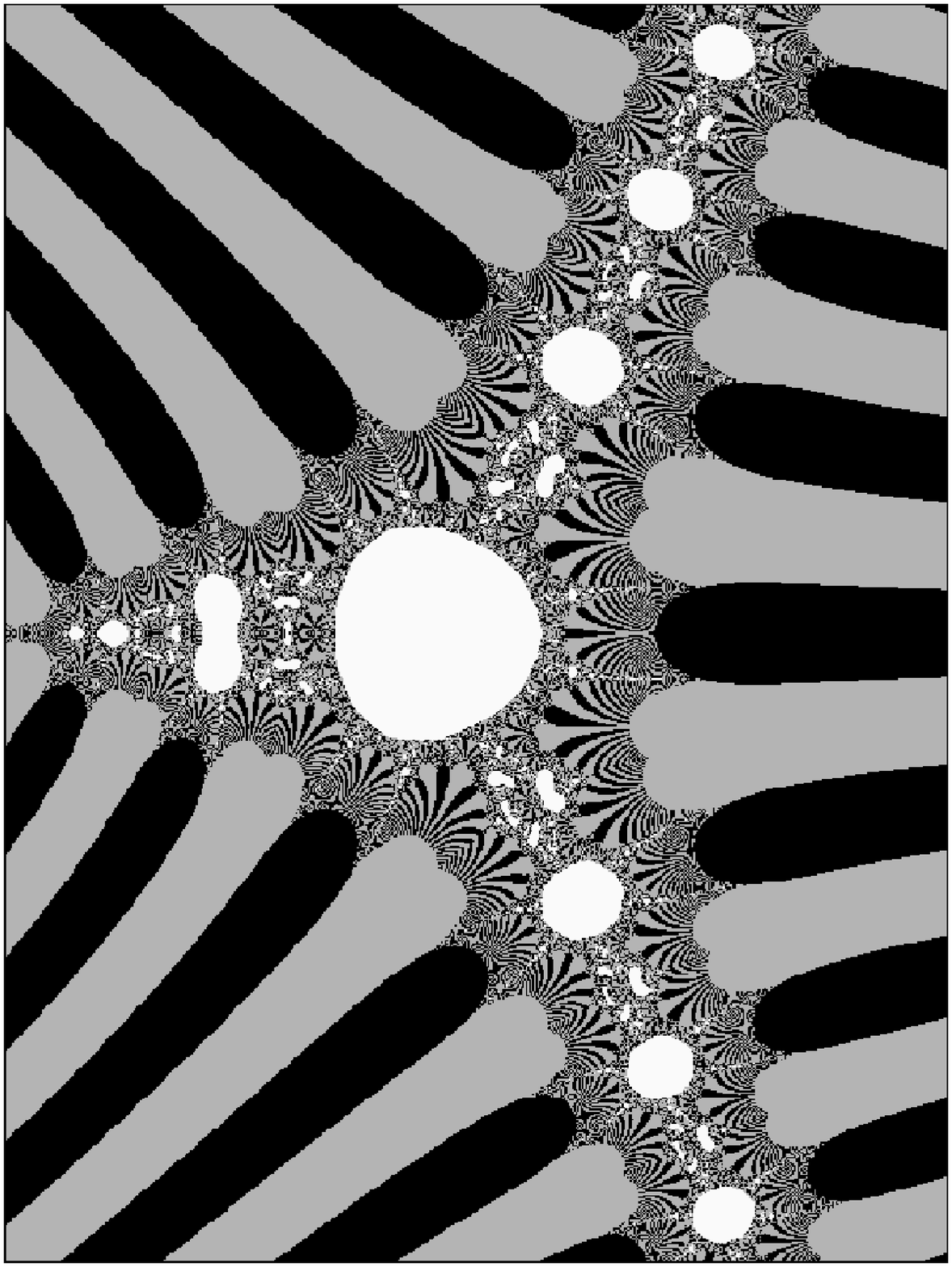}
\quad
\includegraphics[width=4.5cm]{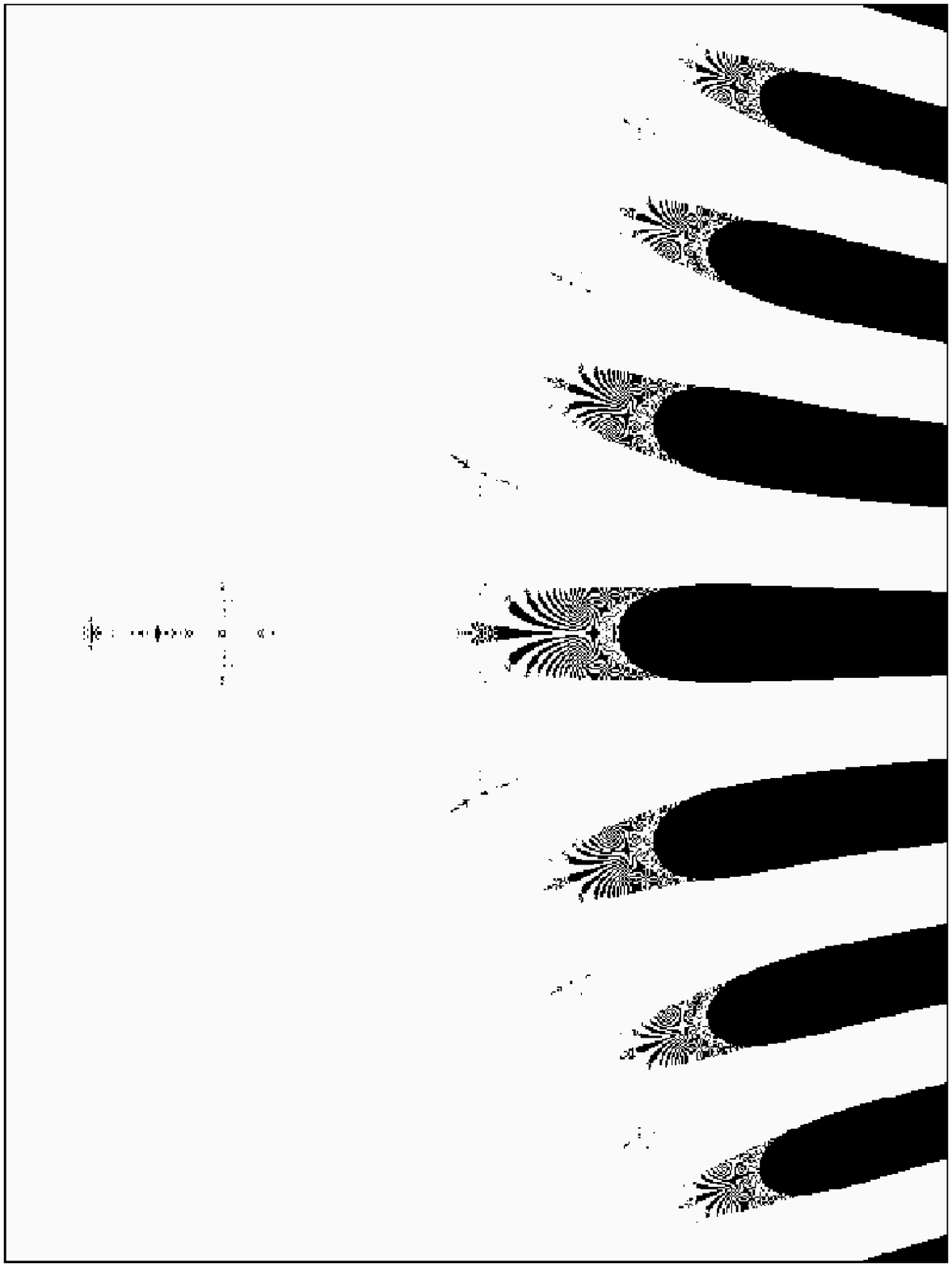}
\quad
\includegraphics[width=4.5cm]{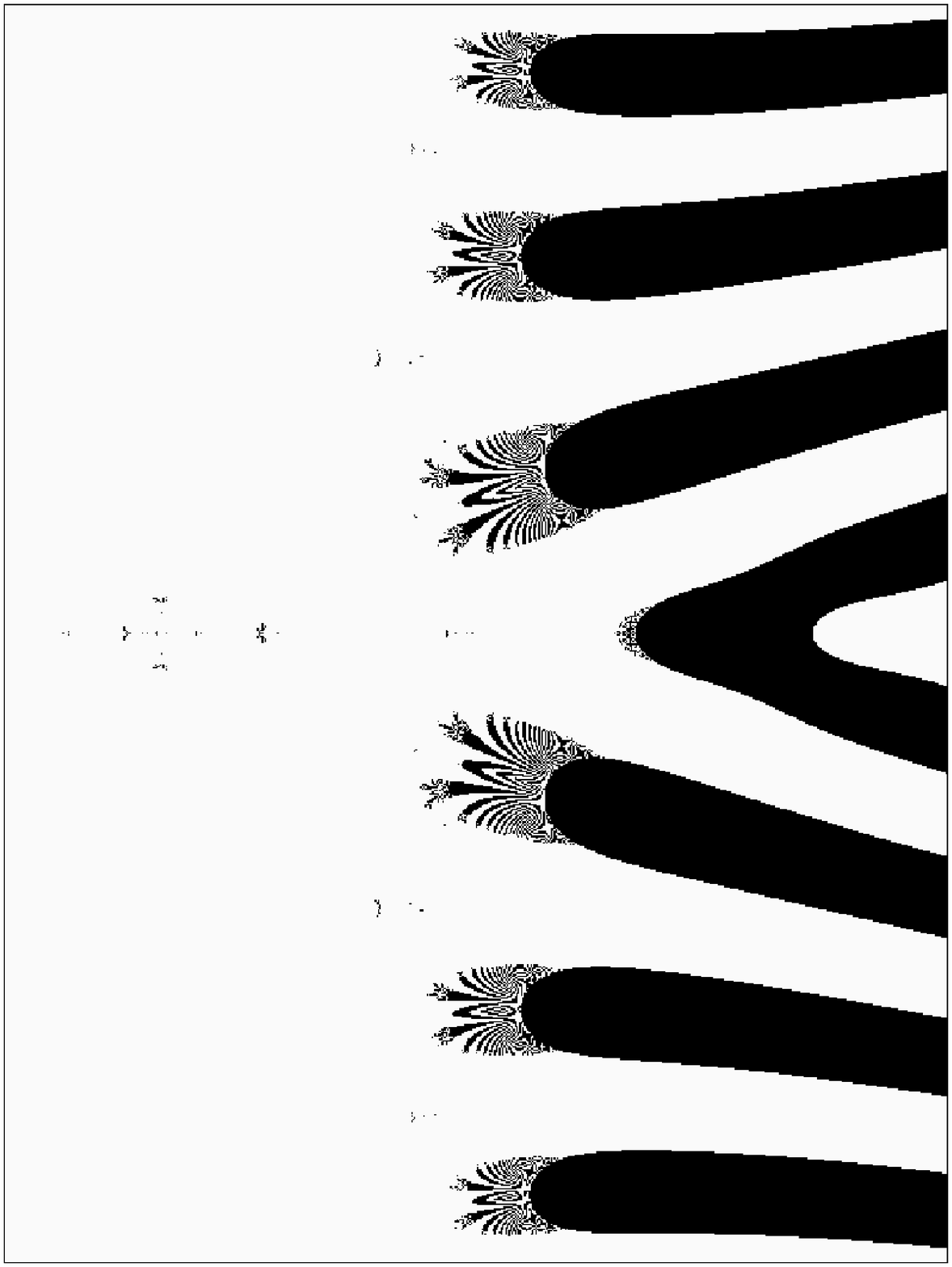}
\caption{The Julia sets of $\Gamma(z)$, $\Gamma(z+1)$ and $\Gamma(z+1)\cos(z)$.}
\label{fig2}
\end{center}
\end{figure}
\end{ex}
\begin{ex} \label{ex3}
Stirling's formula yields that the function 
$f(z)=\Gamma(z+1)\cos(z)$ is bounded on the imaginary
axis and that it is  unbounded in the  right half-plane.
Thus $f$ has a direct tract.
We remark that $f$ satisfies the hypotheses of Proposition~\ref{selberg}.

We note that the tract of $f$ is not logarithmic.
This can be shown directly, but can also be deduced from general
principles.
For example, this follows since a hypothetical logarithmic
tract would be simply connected. By symmetry and since the zeros of 
the cosine are not in a tract, the function $f$ would then have two 
logarithmic tracts in the right half-plane. On the other hand, $f$ 
has order one and hence it follows 
from the Denjoy-Carleman-Ahlfors Theorem~\cite[Section~XI.4]{Nev53}
that $f$ can have at most one direct tract in a half-plane.

The Julia and Fatou sets of $f$ are shown in the right
picture of Figure~\ref{fig2}.
The function has the attracting fixed point $0.6965\ldots$, whose
attracting basin is drawn white, while
the points tending to infinity in the 
direct tract are black.
The range shown is given by $-5\leq \re z\leq 10$, $|\im z|\leq 10$.
\end{ex}

\section{Proof of Theorem~\ref{growth}}
\label{proofgrowth}
We denote the Riesz measure of a subharmonic function $u$
by $\mu_u$; see~\cite[Section~3.5]{HK1}.
\begin{la} \label{lemma1}
Let $u:\C\to [-\infty,\infty)$
be subharmonic and let $\gamma$ be a closed, 
piecewise smooth Jordan curve. 
Denote the interior of $\gamma$ by $\interior (\gamma)$.
Suppose that $\gamma$ has a neighbourhood $U$ where
$u$ has the form $u=\log |f|$ for some 
nonvanishing holomorphic function
$f:U\to \C$. Then
$\mu_u(\interior (\gamma))$ is a nonnegative integer.
\end{la} 
\begin{proof} 
For $z\in U\cup \interior (\gamma)$ we have (see~\cite[Theorem~3.9]{HK1})
$$u(z)=\int_E\log|z-\zeta|d\mu_u(\zeta) +h(z),$$
where $E=\supp(\mu_u)\cap \interior (\gamma)$ is 
compact and $h$ is harmonic in $U\cup \interior (\gamma)$.
We may assume that $U\cup \interior (\gamma)$ is simply connected.
Differentiating we find that if $z\in U$,  then
$$\frac{f'(z)}{f(z)}=\int_E \frac{1}{z-\zeta}d\mu_u(\zeta) +g(z)$$
where $g$ is holomorphic in $U\cup \interior (\gamma)$.
Integrating along $\gamma$ yields
$$\int_\gamma \frac{f'(z)}{f(z)}dz
=\int_E 
\left(\int_\gamma \frac{dz}{z-\zeta}\right)d\mu_u(\zeta)
=2\pi i \int_E d\mu_u(\zeta)
=2\pi i \mu_u(E)
=2\pi i \mu_u(\interior (\gamma)).$$
The integral on the left is purely imaginary, and its imaginary
part is the increase of the argument of $\log f(z)$ as $z$ traverses
the curve $\gamma$. Clearly, this increase is a multiple
of $2\pi$.
The conclusion follows.
\end{proof} 
An immediate consequence of Lemma~\ref{lemma1} is the following result.
\begin{la} \label{lemma1a}
Let $D$ be a direct tract of $f$ and let $v$ be defined
by \textup{(\ref{definitionv})}. Let $K$ be a bounded component 
of the complement of~$D$. Then
$\mu_v(K)\geq 1$.
\end{la} 
\begin{proof} 
It follows from Lemma~\ref{lemma1} that $\mu_v(K)$ is a nonnegative
integer, and the case $\mu_v(K)=0$ is excluded since $v$ is not
harmonic on the boundary of~$K$.
\end{proof} 
The following lemma is essentially Jensen's formula for subharmonic
functions; see~\cite[Section~3.9]{HK1}.
\begin{la} \label{lemma2}
Let $u:\C\to [-\infty,\infty)$
be subharmonic, $a\in\C$ and $r>0$. For $t>0$ put
$n(a,t,u)=\mu_u(\overline{D}(a,t))$. Then
$$\frac{1}{2\pi}\int_0^{2\pi}u(a+re^{i\varphi})d\varphi
=\int_0^r \frac{n(a,t,u)}{t}dt +u(a).$$
\end{la} 
In particular, taking $a=0$ we obtain
a lower bound for $B(r,u)$ from Lemma~\ref{lemma2}.
Another lower bound is given by the following result, which 
in slightly different form
can be found
in~\cite[p.~548]{Hay89} or~\cite[p.~117]{Tsu}.
To state it we  fix a domain $D$ and, for
a circle $C(a,r)=\{z\in\C:|z-a|=r\}$ which
intersects $D$, we denote by $r\theta(a,r)$ the linear measure 
of the intersection. Let $\theta^*(a,r)=\theta(a,r)$ if $C(a,r)\not \subset D$ 
and let $\theta^*(a,r)=\infty$, and thus $1/\theta^*(a,r)=0$, if 
$C(a,r)\subset D$. In this lemma, we use the more general definition 
\[B(r,u)=\max_{|z|=r,z\in D}u(z),\]
since we do not assume that $u$ is defined in~$\C$. 
\begin{la} \label{tsuji1}
Let $D\subset \C$ be an unbounded domain
and suppose that $u:\overline{D}\to [-\infty,\infty)$
is continuous in $\overline{D}$ and subharmonic in~$D$.
Suppose also that $u$ is bounded above
on $\partial D$, but not bounded above in~$D$.
Let $0<\kappa<1$
and let $r_0>0$ be such that $C(0,r_0)$ intersects~$D$.
Then 
$$\log
B(r,u)
\geq\pi\int_{r_0}^{\kappa r}\frac{dt}{t\theta^*(0,t)}-\OO(1)$$
as $r\to\infty$.
\end{la} 
\begin{proof}[Proof of Theorem~\ref{growth}]
As already mentioned after the statement of Theorem~\ref{growth},
(\ref{growtha}) follows from~(\ref{growthB}).
Thus it suffices to prove~(\ref{growthB}).

Suppose first that
the complement of $D$ has infinitely many bounded components.
It then follows from Lemma~\ref{lemma1a} 
that the function $n(0,t,v)$ defined in
Lemma~\ref{lemma2} tends to $\infty$ as $t$ tends to $\infty$.
By Lemma~\ref{lemma2}  we have
$$B(r,v)\geq\frac{1}{2\pi}\int_0^{2\pi}v(re^{i\varphi})d\varphi
=\int_0^r \frac{n(0,t,v)}{t}dt +v(0)\geq
\int_{\sqrt{r}}^r \frac{n(0,t,v)}{t}dt
\geq\tfrac12 n(0,\sqrt{r},v)\log r,$$
and~(\ref{growthB}) follows.

Suppose now that
the complement of $D$ has only finitely many bounded components.
Since 
the complement of $D$ is unbounded 
this implies that the complement of $D$ has at least one 
unbounded component.
For the function $\theta^*$ defined before Lemma~\ref{tsuji1}
we thus have $\theta^*(0,t)\leq 2\pi$ for all large $t$, say for 
$t\geq t_0>1$. 
Lemma~\ref{tsuji1} now yields
$$\log
B(r,v)
\geq\pi\int_{t_0}^{\kappa r}\frac{dt}{t\theta^*(0,t)}-\OO(1)
\geq \tfrac12 \log r -\OO(1)$$
and hence 
$$B(r,v)\geq c\sqrt{r}$$
for some $c>0$ and all large~$r$.
Again~(\ref{growthB}) follows.
\end{proof} 
\section{Some estimates of harmonic measure}
\label{harmonicmeasure}
In this section and the next we give some auxiliary results that
are needed in the proof of Theorem~\ref{th1}.
The following estimate of harmonic measure can be found 
in~\cite[p.~112]{Tsu}. 
We mention that
Lemma~\ref{tsuji1} is deduced from this estimate 
in~\cite{Tsu}. 
We use the standard notation for harmonic measure 
and the terminology introduced before Lemma~\ref{tsuji1}.
\begin{la} \label{tsuji2}
Let $D\subset \C$ be a domain, $a\in D$ and $r>0$. 
Let $V$ be the component of $D\cap D(a,r)$ 
that contains $a$ and let $\Gamma=\partial V\cap C(a,r)$.
For  $0<\kappa<1$ we then have
$$\omega(a,\Gamma,V)\leq\frac{3}{\sqrt{1-\kappa}}
\exp\left(-\pi \int_0^{\kappa r}\frac{dt}{t\theta^*(a,t)}\right).$$
\end{la} 
The following result is known as the two constants theorem; 
see~\cite[Section~III.2]{Nev53}, although our version is somewhat
different from the one given there.
\begin{la} \label{twoconstants}
Let $V$ be a bounded domain with piecewise smooth boundary.
Let $\Sigma$ be a subset of $\partial V$ consisting of
finitely many boundary arcs and let $m,M$ be real constants
with $m<M$.
Suppose that $u:\overline{V}\to [-\infty,\infty)$
is continuous in $\overline{V}$ and subharmonic in~$V$.
Suppose also that $u(z)\leq M$ for all $z\in \overline{V}$
and that  $u(z)\leq m$ for $z\in \Sigma$.
Then
$$u(z)\leq \omega(z,\Sigma,V)m +(1-\omega(z,\Sigma,V))M.$$
\end{la} 
To prove the lemma one has only to note that the right-hand
side is a harmonic function which majorizes $u$ on the boundary
of~$V$. Since $u$ is subharmonic
it thus majorizes $u$ also in the interior of~$V$.
\section{Growth lemmas for real functions}
The following result is a version of a classical 
lemma due to Borel and Nevanlinna; see~\cite{Nev31} or~\cite[Section~3.3]{ChYe}.
\begin{la} \label{growthlemma1}
Let $x_0>0$ and let $T:
[x_0,\infty)\to (0,\infty)$ be nondecreasing.
Let $0<\alpha<\beta$. 
Then there exists a set $E\subset [x_0,\infty)$ of finite measure
such that if $x\notin E$, then
\begin{equation}
T\left(x+T(x)^{-\beta}\right) <  
\left(1+T(x)^{-\alpha}\right) T(x)     \label{upperboundT}
\end{equation}
and
\begin{equation}
T\left(x-T(x)^{-\beta}\right) >     
\left(1-T(x)^{-\alpha}\right) T(x).     \label{lowerboundT}
\end{equation}
\end{la}
\begin{proof}
Let $E_1$ be the subset of $[x_0,\infty)$ where~\eqref{upperboundT}
fails and let $E_2$ be the one  where~\eqref{lowerboundT} fails,
so that we can take $E=E_1\cup E_2$.
To estimate the size of $E_1$ we may assume that $E_1$ is unbounded.
We choose $x_1\in E_1\cap [\inf E_1,\inf E_1 +\tfrac12]$ and 
put $x_1'=x_1+T(x_1)^{-\beta}$.
Recursively we then choose
$$x_j\in E_1\cap \left[\inf \left(E_1\cap [x_{j-1}',\infty)\right),
\inf \left(E_1\cap [x_{j-1}',\infty)\right)+2^{-j}\right]$$
and put $x_j'=x_j+T(x_j)^{-\beta}$.
Then 
$$T(x_{j+1})\geq T(x_{j}')=T\left(x_j+T(x_j)^{-\beta}\right)
\geq \left(1+T(x_j)^{-\alpha}\right) T(x_j).$$
It follows from this that $T(x_j)\to\infty$, since otherwise 
there exists $M>0$ with $T(x_j)\leq M$ for all $j\in\N$ and thus
$$M\geq T(x_{j+1})\geq\left(1+M^{-\alpha}\right) T(x_j)
\geq\left(1+M^{-\alpha}\right)^j T(x_1),$$
for all $j\in\N$, a contradiction.

We shall prove by induction that
there exists a constant $c_1>0$ such that 
\begin{equation} \label{growthT}
T(x_j)\geq c_1 j^{1/\alpha}
\end{equation}
for all $j\in\N$.

In order to do so we note first 
\begin{equation} \label{conditionforc}
\left(\frac{j+1}{j}\right)^{1/\alpha}
=\left(1+\frac{1}{j}\right)^{1/\alpha}
\leq 1+\frac{c_1^{-\alpha}}{j}
\end{equation}
for all $j\in\N$, if $c_1$ is chosen small enough, say for $c_1\leq c'$.
Next we note that there exists $t_\alpha>0$ such
that the function $t\mapsto (1+t^{-\alpha})t$ increases
for $t\geq t_\alpha$. 

We now choose ${j_0}$ so large that $T(x_{j_0})\geq t_\alpha$
and $c'{j_0}^{1/\alpha}\geq t_\alpha$.
Choosing $c_1$ such that $c_1{j_0}^{1/\alpha}=t_\alpha$
we see that~\eqref{growthT} holds for $j=j_0$. 
Suppose now 
that~\eqref{growthT} holds for some $j\geq j_0$.
Then 
$$T(x_{j+1})\geq\left(1+T(x_j)^{-\alpha}\right) T(x_j)
\geq\left(1+\left(c_1 j^{1/\alpha}\right)^{-\alpha}\right)c_1 j^{1/\alpha}
=\left( 1+\frac{c_1^{-\alpha}}{j}\right)c_1 j^{1/\alpha}.$$
Combining this with~\eqref{conditionforc} we obtain 
$$T(x_{j+1})\geq c_1 (j+1)^{1/\alpha},$$
and thus we see that~\eqref{growthT} holds for all $j\geq j_0$.
Adjusting the value of $c_1$ if necessary we
may thus assume that~\eqref{growthT} is satisfied for all $j\in\N$.

Since $T(x_j)\to\infty$ and thus $x_j\to\infty$ as $j\to\infty$ we have
$$E_1\subset \bigcup_{j=1}^\infty \left[x_j-2^{-j},x_j'\right]$$
and hence 
$$\meas E_1\leq \sum_{j=1}^\infty \left(x_j'-x_j+2^{-j}\right)
=\sum_{j=1}^\infty\frac{1}{T(x_j)^\beta}
+\sum_{j=1}^\infty \frac{1}{2^{j}}
\leq c_1^{-\beta}\sum_{j=1}^\infty\frac{1}{j^{\beta/\alpha}}+1
<\infty.$$
To estimate $E_2$ we proceed similarly. We may assume that 
$E_2\neq \emptyset$ and fix $R>x_0$ so large that $E_2\cap [x_0,R]
\neq \emptyset$.
We choose
$$z_1\in
E_2\cap \left[\sup \left(E_2\cap [x_0,R]\right)-\tfrac12,
\sup \left(E_2\cap [x_0,R]\right)\right]$$
and put
$z_1'=z_1-T(z_1)^{-\beta}$.
Recursively we then choose
$$z_j\in E_2\cap \left[\sup \left(E_2\cap [x_0,z_{j-1}']\right)-2^{-j},
\sup \left(E_2\cap [x_0,z_{j-1}']\right)\right]$$ 
and
put $z_j'=z_j-T(z_j)^{-\beta}$,
as long as $E_2\cap [x_0,z_{j-1}'] \neq \emptyset$.
However, since 
\begin{eqnarray*}
T(z_{j+1})
&\leq& 
T(z_{j}')\\
&=&
T\left(z_j-T(z_j)^{-\beta}\right)\\
&\leq&
\left(1-T(z_j)^{-\alpha}\right) T(z_j)\\
&\leq&
\left(1-T(z_1)^{-\alpha}\right) T(z_j)\\
&\leq&
\left(1-T(z_1)^{-\alpha}\right)^j T(z_1)
\end{eqnarray*}
we see that the process stops and we obtain two finite
sequences $(z_1,\dots,z_N)$ and $(z_1',\dots,z_N')$ with 
$$E_2\cap [x_0,R] \subset \bigcup_{j=1}^N [z_j',z_j+2^{-j}].$$
With $y_j=z_{N-j+1}$ we thus have
$$E_2\cap [x_0,R] \subset \bigcup_{j=1}^N [y_j',y_j+2^{j-N-1}]$$
and
$$T(y_{j})\leq \left(1-T(y_{j+1})^{-\alpha}\right) T(y_{j+1}).$$
As in the proof of~\eqref{growthT}
we see that for sufficiently small $c_2$  and large~$j$ the conditions 
$T(y_j)\geq c_2 j^{1/\alpha}$ and $T(y_{j+1})< c_2 (j+1)^{1/\alpha}$
are incompatible for large $j$,
and as before we can deduce that there exists a positive constant~$c_2$ 
independent of $R$ such that 
$$T(y_j)\geq c_2 j^{1/\alpha}$$
for $1\leq j\leq N$.
We conclude that
$$\meas \left(E_2\cap[x_0,R]\right)\leq\sum_{j=1}^N (y_j-y_j'+2^{j-N-1})
=\sum_{j=1}^N \frac{1}{T(y_j)^\beta} +\sum_{j=1}^N 2^{j-N-1}
\leq c_2^{-\beta}\sum_{j=1}^\infty\frac{1}{j^{\beta/\alpha}}+1,$$
and thus
$\meas E_2 <\infty.$
\end{proof}
We will apply this lemma to the derivative of a convex function~$\Phi$.
We note here that a convex function has one-sided derivatives 
at all points and is differentiable except for at most countably many points. 
In the following result it will be irrelevant how we define the 
derivative at these countably many points, but to be definite we denote
by $\Phi'$ the derivative from the right, which is then an increasing
function.
\begin{la} \label{growthlemma2}
Let $x_0>0$ and let $\Phi:
[x_0,\infty)\to (0,\infty)$ be increasing and convex.
Let $\beta>\frac12$.
Then there exists a set $E\subset [x_0,\infty)$ of finite measure
such that
\begin{equation} \label{2e}
\Phi(x+h)\leq\Phi(x)+\Phi'(x)h+o(1) \quad \text{for}\quad 
|h|\leq \Phi'(x)^{-\beta},
\quad x\notin E,
\end{equation}
uniformly as $x\to\infty$.
\end{la}
For differentiable $\Phi$ and for $h$ in the range
$|h|\leq C/\Phi'(x)$ with a constant $C>0$ this result was 
proved in~\cite[Lemma 2]{Ber90}.
\begin{proof}[Proof of Lemma~\ref{growthlemma2}]
Since $\Phi'$ is nondecreasing, $\lim_{x\to\infty} \Phi'(x)$
exists. 
It is easy to see that~\eqref{2e} holds for all large $x$ if
this limit is finite. Hence we assume that 
$\lim_{x\to\infty} \Phi'(x)=\infty$. 

We apply Lemma \ref{growthlemma1} with $T=\Phi'$ and some $\alpha$
satisfying $\frac12<\alpha<\beta$.
For $x\notin E$ and $0<h\leq \Phi'(x)^{-\beta}$ we then have
\begin{eqnarray*}
\Phi(x+h)
&=&
\Phi(x)+\int_x^{x+h} \Phi'(u)du\\
&\leq &
\Phi(x)+\Phi'(x+h)h\\
&\leq &
\Phi(x)+\Phi'(x+\Phi'(x)^{-\beta})h\\
&\leq& 
\Phi(x)+\left(1+\Phi'(x)^{-\alpha}\right) \Phi'(x) h\\
&=&
\Phi(x)+\Phi'(x) h+\Phi'(x)^{1-\alpha}h\\
&\leq &
\Phi(x)+\Phi'(x) h+\Phi'(x)^{1-\alpha-\beta}\\
&\leq&
\Phi(x)+\Phi'(x) h+o(1)
\end{eqnarray*}
as $x\to\infty$.
The case $-\Phi'(x)^{-\beta}\leq h<0$ is analogous.
\end{proof}
We apply Lemma~\ref{growthlemma2}
to $\Phi(x)=B(e^x,v)$ where $v$ is subharmonic
and $B$ is defined by~\eqref{defBrv}.
Then $\Phi'(x)=e^xB'(e^x,v)=a(e^x,v)$.
With $r=e^x$ and $s=e^{x+h}$ we obtain
$$B(s,v)=\Phi(x+h)\leq\Phi(x)+\Phi'(x)h+o(1)
=B(r,v)+a(r,v)\log\frac{s}{r}+o(1)$$
for $r\notin F=\exp E$, provided
$|\log(s/r)|=|h|\leq \Phi'(x)^{-\beta}=a(r,v)^{-\beta}$.

Hence we obtain the following result.
\begin{la} \label{growthlemma3}
Let $v:\C\to [-\infty,\infty)$
be subharmonic and let $\beta>\tfrac12$.
Then there exists a set $F\subset [1,\infty)$ of finite logarithmic 
measure such that 
$$B(s,v)\leq B(r,v)+a(r,v)\log\frac{s}{r}+o(1)
\quad \text{for}\quad 
\left|\log\frac{s}{r}\right|\leq \frac{1}{a(r,v)^\beta},
\quad r\notin F,$$
uniformly as $r\to\infty$.
\end{la}

\section{Proof of Theorem~\ref{th1}}
\label{proof1}
We apply Lemma~\ref{growthlemma3} for some $\beta<1$ satisfying 
$\tfrac12<\beta<\tau$ and we apply
Lemma~\ref{growthlemma5} for some $\eps>0$
such that $(1-\beta)(1+\eps)<1$.
Let $F$ be the union of the exceptional sets of these lemmas.
We put $\rho=2 r a(r,v)^{-\tau}$ whenever $r$ is so large that $a(r,v)\neq 0$.

We consider the function
\begin{equation}\label{definitionu}
u(z)=v(z)-v(z_r)-a(r,v)\log \frac{|z|}{r}
=v(z)-B(r,v) -a(r,v)\log \frac{|z|}{r}.
\end{equation}
For $z\in \overline{D}(z_r,512 \rho)$ we have 
$$\left| \frac{z-z_r}{z_r}\right| \leq
\frac{512 \rho}{r}=\frac{1024}{a(r,v)^{\tau}}= o(1)$$
as $r\to\infty$ by~(\ref{growtha}) and thus
\begin{equation}\label{logzr}
\left| \log \frac{|z|}{r}\right|
=\left| \log \left|1+\frac{z-z_r}{z_r}\right|\right|
\leq 2\left| \frac{z-z_r}{z_r}\right|
\leq \frac{2048}{a(r,v)^{\tau}}
\leq \frac{1}{a(r,v)^{\beta}}
\end{equation}
for large~$r$.  Since
$$u(z)\leq B(|z|,v)-B(r,v)-a(r,v)\log \frac{|z|}{r}$$
by the definition of $u$ 
we conclude from Lemma~\ref{growthlemma3} that 
\begin{equation}\label{boundu}
u(z)\leq o(1)
\quad \text{for} \quad z\in \overline{D}(z_r,512 \rho), \quad r\notin F,
\end{equation}
as $r\to\infty$.

We now show that $D(z_r, \rho)\subset D$ if $r\notin F$ is sufficiently large.
Suppose that there exists $\xi\in D(z_r, \rho)$ with $\xi\notin D$.
We denote by $K$ the component of the complement of $D$ that contains $\xi$
and distinguish two cases.

{\em Case 1.} $K\not\subset D(z_r,256 \rho)$.

Then $K$ intersects  $\partial D(z_r,t)$ for
$\rho\leq t\leq 256 \rho$. 
Thus $\theta^*(z_r,t)\leq 2\pi$ for $\rho\leq t\leq 256 \rho$.
Let $V$ be the component of 
$D\cap  D(z_r,512\rho)$ that contains $z_r$ and let 
$\Gamma=\partial V\cap \partial D(z_r,512\rho)$.
Lemma~\ref{tsuji2} yields
with  $\kappa=\frac12$ that
$$\omega(z_r,\Gamma,V)\leq
3\sqrt{2}\exp\left(-\pi \int_{\rho}^{256 \rho}\frac{dt}{t\theta^*(z_r,t)}\right)
\leq 3\sqrt{2} \exp\left(-\frac12 \int_{\rho}^{256 \rho} \frac{dt}{t}\right)
=\frac{3\sqrt{2}}{2^4}<\frac12.$$
For $\Sigma=\partial V\setminus \Gamma$ we thus have
\begin{equation}\label{omegaSigma}
\omega(z_r,\Sigma,V)=1-\omega(z_r,\Gamma,V)> \frac12.
\end{equation}
For $z\in\Sigma$ we have $v(z)=0$ and thus
we deduce from~\eqref{definitionu},~\eqref{logzr} and~\eqref{aB}
that if $r\notin F$ is sufficiently large, then
\begin{eqnarray*}
u(z)
&=&
-B(r,v)-a(r,v)\log \frac{|z|}{r}\\
&\leq& -B(r,v)+a(r,v)^{1-\beta}\\
&\leq& -B(r,v)+B(r,v)^{(1-\beta)(1+\eps)}\\
&\leq& -\tfrac12 B(r,v).
\end{eqnarray*}
By this inequality and~\eqref{boundu},
we can apply Lemma~\ref{twoconstants},
for large $r\notin F$,
with $m=-\tfrac12 B(r,v)$ and $M=1$,
and we obtain, by~\eqref{omegaSigma},
$$u(z_r)\leq -\tfrac12\omega (z_r,\Sigma,V) B(r,v)+1-\omega (z_r,\Sigma,V)
\leq -\tfrac14 B(r,v) + 1.$$
This is a contradiction since
$u(z_r)=0$ by the definition of $u$,
but $B(r,v)\to\infty$ as $r\to\infty$.

{\em Case 2.} $K\subset D(z_r,256 \rho)$.

We again use the function $u$ defined by~(\ref{definitionu}).
For large $r$ we have $0\notin D(z_r,512 \rho)$ so that 
the difference of $u$ and $v$ is harmonic
in $D(z_r,512 \rho)$. Hence their Riesz measures
in this disc coincide.
Lemma~\ref{lemma1a} implies that the 
Riesz measure of $K$ is a positive integer
and thus $n(z_r,t,u)=n(z_r,t,v)\geq 1$ for $t\geq 256 \rho$.
Hence
$$\int_0^{512 \rho} \frac{n(z_r,t,u)}{t}dt \geq 
\int_{256 \rho}^{512 \rho} \frac{dt}{t}
=\log 2.$$
On the other hand,
Lemma~\ref{lemma2}  and~(\ref{boundu}) yield
$$\int_0^{512 \rho} \frac{n(z_r,t,u)}{t}dt
= 
\frac{1}{2\pi}\int_0^{2\pi}u(z_r+512 \rho e^{i\varphi})d\varphi
-u(z_r)\leq o(1)$$
as $r\to\infty$,  $r\notin F$.
The last two inequalities yield a contradiction.

This completes the proof that $D(z_r,\rho)\subset D$ 
for large $r\notin F$. 
In order to prove~(\ref{wvlike}) we note that
since $D(z_r,\rho)\subset D$ we may define 
a holomorphic function $g:D(z_r,\rho)\to\C$ by
$$ g(z)=\log\frac{f(z)}{f(z_r)}-a(r,v)\log\frac{z}{z_r} 
=\log \left( \frac{f(z)}{f(z_r)}  
\left(\frac{z_r}{z} \right)^{a(r,v)} \right),$$
with the branches of the
logarithms chosen such that $g(z_r)=0$.
By the Borel-Carath\'eodory inequality (see, e.g.,~\cite[p.~20]{Val23}),
we have 
$$\max_{|z-z_r|=t}|g(z)|\leq 4\max_{|z-z_r|=2t}\re g(z)$$
for $0<t<\rho/2$.
Since $\re g(z)=u(z)$ we can now deduce from~(\ref{boundu}) that 
if $z\in D(z_r, \rho/2)=D(z_r, r a(r,v)^{-\tau})$,
then $g(z)\to 0$ as $r\to\infty$, $r\notin F$.
This yields~(\ref{wvlike}).

It follows from~(\ref{wvlike}) and Lemma~\ref{growthlemma3} that,
as $r\to\infty$, $r\notin F$, we have
$$M_D(|z|)\geq |f(z)|\geq (1-o(1)) 
\left|\frac{z}{z_r}\right|^{a(r,v)}|f(z_r)|$$
and
$$M_D(|z|)=\exp B(|z|,v)\leq \exp\left(B(r,v)+a(r,v)\log\frac{|z|}{r}+o(1)\right)
=(1+o(1))\left|\frac{z}{z_r}\right|^{a(r,v)}|f(z_r)|$$
for $z\in D(z_r, r a(r,v)^{-\tau})$.
The last two inequalities give~(\ref{wvlike4}).  

\section{Applications to differential equations}
\label{diffeq}
Wiman-Valiron theory has found many applications in the
theory of differential equations in the complex domain,
see, e.g.,~\cite{JV,Lai,Wit}. 
It seems plausible that~\eqref{wvlike3} allows us to extend some of these results 
to meromorphic solutions with a direct singularity.
We confine ourselves to giving one example of an application
to differential equations.
Another application of Theorem~\ref{th1} to complex differential equations is
given in~\cite{ELN}.
In order to state our application, let $f$ be meromorphic,
$n\in\N$, $t_j\in \N_0 =\N\cup\{0\}$
for $j=0,1,2,\dots,n$, and put $t=(t_0,t_1,\dots,t_n)$.
Define $M_t[f]$ by
\[
M_t[f](z)=f(z)^{t_0}f'(z)^{t_1}f''(z)^{t_2}\dots f^{(n)}(z)^{t_n},
\]
with the convention that $M_{(0)}[f]=1$. We call
$d(t)=t_0+t_1+\dots+t_n$ the {\em degree} and 
$w(t)=t_1+2t_2+\dots+nt_n$ the {\em weight} of $M_t[f]$.
An {\em algebraic differential equation} is an equation of the form
\begin{equation} \label{ade}
\sum_{t\in T}c_tM_t[f]=0,
\end{equation}
where the $c_t$ are polynomials and $T$ is a finite
index set. 
\begin{thm} \label{th5}
Let $f$ be a transcendental 
meromorphic solution of \textup{(\ref{ade})}
with a direct singularity over $\infty$.
Let $S$ be the set of all $s\in T$ for which $d(s)=\max_{t\in T} d(t)$.
Then $S$ has at least two elements.

Let $\Lambda$ be the set of all $\lambda\in\N$ for which there
exists $s\in S$ satisfying $w(s)=\lambda$.
Suppose that
$$\sum_{t\in  S, w(t)=\lambda}c_t\neq 0$$
for all $\lambda\in\Lambda$. 
Then the order of $f$ is at least $1/\max(\Lambda)$.
\end{thm}
For an entire function $f$ this result can be 
found in~\cite[p.~64, p.~71]{Wit}.
The proof based on Wiman-Valiron theory extends to our 
setting without difficulty, but for completeness we include 
the argument below.

For an entire solution $f$, Wiman-Valiron theory also yields an upper
bound for the growth of~$f$.
In particular, an entire solution has finite order. 
For a meromorphic solution with a direct singularity over $\infty$
we still find that the order in a direct tract is bounded, 
but the method does not give any information about what happens
outside the tract.

\begin{proof}[Proof of Theorem~\ref{th5}]
Let $D$ be a direct tract of $f$, let $R$ be as in Definition~\ref{defin1}
and define $v$ by~(\ref{definitionv}).
We use~(\ref{wvlike3}) with $z=z_r$ and obtain
\begin{equation}\label{ade1}
\sum_{t\in T}
c_t(z_r) 
f(z_r)^{d(t)}
\left(\frac{a(r,v)}{z_r}\right)^{w(t)}
\left(1+\varepsilon_t(r)\right)
=0,
\end{equation}
where $\varepsilon_t(r)\to 0$ as $r\to\infty$, $r\notin F$.
By Lemma~\ref{growthlemma5} we may assume that 
$a(r,v)\leq B(r,v)^2=\left(\log |f(z_r)|\right)^2$ for $r\notin F$
and this implies that as $r\to\infty$, $r\notin F$,
$$
c_t(z_r) 
f(z_r)^{d(t)}
\left(\frac{a(r,v)}{z_r}\right)^{w(t)}
=
o\left(
c_s(z_r) 
f(z_r)^{d(s)}
\left(\frac{a(r,v)}{z_r}\right)^{w(s)}
\right)
$$
if $d(t)<d(s)$. 
After dividing by $f(z_r)^d$ where $d=\max_{t\in T} d(t)$ the equation
(\ref{ade1}) thus takes the form
\begin{equation}\label{ade2}
\sum_{t\in S}
c_t(z_r) 
\left(\frac{a(r,v)}{z_r}\right)^{w(t)}
\left(1+\varepsilon^*_t(r)\right)
=0,
\end{equation}
where $\varepsilon^*_t(r)\to 0$ as 
$r\to\infty$, $r\notin F$.
This implies that $S$ has at least two elements.

By hypothesis,
$$u_\lambda=\sum_{s\in  S, w(s)=\lambda}c_s\neq 0,$$
and hence there exist $b_\lambda\in\C\setminus\{0\}$ and
$d_\lambda\in\N_0$ such that $u_\lambda(z)\sim b_\lambda
z^{d_\lambda}$ as $z\to\infty$.
Equation~(\ref{ade2}) now takes the form
$$
\sum_{\lambda\in \Lambda}
b_\lambda
z_r^{d_\lambda-\lambda}
a(r,v)^\lambda
\left(1+\varepsilon^{**}_t(r)\right)
=0,
$$
where $\varepsilon^{**}_t(r)\to 0$ as
$r\to\infty$, $r\notin F$. Except for the terms $\varepsilon^{**}_t(r)$
and the exceptional set $F$, this is an algebraic equation for $a(r,v)$.
This makes it plausible that $a(r,f)$ grows like a solution of 
the associated algebraic equation. 
It turns out that this is indeed the case. 
We omit the argument justifying this since it can be found in detail
in~\cite[Hilfssatz 22.2]{JV}.
It follows from the argument given there 
that there exists a positive real number $\tau$ and a 
positive rational number $\kappa$ such that 
\begin{equation}\label{ade4}
a(r,v)\sim \tau r^\kappa
\end{equation}
as $r\to\infty$.
In fact, the possible values for $\kappa$ can be computed
from the Newton-Puiseux diagram associated to the algebraic 
equation. In particular it follows that $\kappa\geq 1/\max(\Lambda)$.
Integrating~(\ref{ade4}) yields 
$$
B(r,v)\sim \frac{\tau}{\kappa}  r^\kappa.
$$
By~\cite[Theorem 3.19]{HK1} we have
$$
B(\rho,v)
\leq
\left(\frac{r+\rho}{r-\rho}\right)
\frac{1}{2\pi}\int_0^{2\pi} v(r e^{i\varphi}) d\varphi
$$
for $0<\rho<r$.
Now
$$\frac{1}{2\pi}\int_0^{2\pi} v(r e^{i\varphi}) d\varphi
\leq
m\left(r,\frac{f}{R}\right)
\leq T\left(r,\frac{f}{R}\right)
\leq T(r,f)+O(1)$$
as $\rho\to\infty$,
where $m(r,f)$ and $T(r,f)$ are the usual quantities from
Nevanlinna theory. Combining the last three estimates 
and choosing $r=2\rho$
we obtain
$$ T(r,f)\geq \frac13 B\left(\frac{r}{2},v\right)
\geq (1-o(1)) 
\frac{\tau}{3\kappa 2^\kappa}  r^\kappa
$$ 
as $r\to\infty$. Thus the order of $f$ is at least~$\kappa$.
\end{proof}

\end{document}